\documentclass[10pt,twoside,fleqn,a4paper]{article}

\usepackage{amssymb}
\usepackage{epsfig}

\newtheorem{theorem}{Theorem}[section]
\newtheorem{proof}{Theorem}[section]

\title{
  Geometrically Graded $h$-$p$ Quadrature Applied to the
  Complex Boundary Integral Equation Method for the Dirichlet
  Problem with Corner Singularities
}

\author{David De Wit}

\date{November 1992}

\begin{document}

\maketitle

\begin{abstract}
  Boundary integral methods for the solution of boundary value
  PDEs are an alternative to `interior' methods, such as finite
  difference and finite element methods. They are attractive on
  domains with corners, particularly when the solution has
  singularities at these corners. In these cases, interior
  methods can become excessively expensive, as they require a
  finely discretised 2D mesh in the vicinity of corners, whilst
  boundary integral methods typically require a mesh discretised
  in only one dimension, that of arc length.

  Consider the Dirichlet problem. Traditional boundary integral
  methods applied to problems with corner singularities involve a
  (real) boundary integral equation with a kernel containing a
  logarithmic singularity.  This is both tedious to code and
  computationally inefficient.  The CBIEM is different in that it
  involves a complex boundary integral equation with a smooth
  kernel. The boundary integral equation is approximated using a
  collocation technique, and the interior solution is then
  approximated using a discretisation of Cauchy's integral
  formula, combined with singularity subtraction.

  A high order quadrature rule is required for the solution of
  the integral equation. Typical corner singularities are of
  square root type, and a `geometrically graded $h$-$p$'
  composite quadrature rule is used.  This yields efficient,
  high order solution of the integral equation, and thence the
  Dirichlet problem.

  Implementation and experimental results in \textsc{matlab} code
  are presented.
\end{abstract}

\pagebreak

%%%%%%%%%%%%%%%%%%%%%%%%%%%%%%%%%%%%%%%%%%%%%%%%%%%%%%%%%%%%%%%%%%%%%%%%

\section{Introduction}

This report describes a research project carried out from March to
October 1992, at the Department of Mathematics, The University of
Queensland, Australia. It was carried out under the supervision of Dr
Graeme A. Chandler, and was accredited as a \#30 project, coded MN882.

Techniques related to the CBIEM have been analysed
in~\cite{ProssdorfRathsfeld:89}, and used in~\cite{Forbes:85,Dold:92}.
The CBIEM is also closely related to the `Complex Variable Boundary
Element Method'~\cite{Hromadka:84}. This report contains an application
of it, using $h$-$p$ quadrature to achieve high rates of
convergence, even in the presence of corner singularities. This
application owes its conception to my supervisor.

The ideas of graded meshes and $h$-$p$ quadrature (numerical
integration) are presented in \S\ref{sec:quadrature}, and are
illustrated by experimental results. The CBIEM itself is described in
\S\ref{sec:CBIEM}.  \S\ref{sec:results} details numerical
implementation of the CBIEM, using the quadrature technique described
in \S\ref{sec:quadrature}, and presents error results for some test
problems.  \S\ref{sec:further} concludes the report with suggestions
for further development. \textsc{matlab} code written for the
implementation is listed in Appendix \ref{app:code}.

\pagebreak

%%%%%%%%%%%%%%%%%%%%%%%%%%%%%%%%%%%%%%%%%%%%%%%%%%%%%%%%%%%%%%%%%%%%%%%%

\section{$h$-$p$ Quadrature Methods}
\label{sec:quadrature}

\subsection{Introduction}

This section describes a high order numerical integration (quadrature)
technique, that retains its high order in the case of end point
singularities in the integrand. The method uses a \emph{graded mesh},
with integration rules of high order used on larger intervals, and low
order on smaller intervals. To achieve the `best' possible
convergence rates, whilst including the end points of each interval,
the basic quadrature rules used are Gau{\ss}--Lobatto. The underlying
mesh is graded in a \emph{geometric} manner. As the method of using
different quadrature rules on internal intervals is a generalisation of
earlier `$h$' and `$p$' methods, the resulting composite
quadrature rule is called a `geometrically graded $h$-$p$'
method~\cite{BabuskaDorr:81}.

%%%%%%%%%%%%%%%%%%%%%%%%%%%%%%%%%%%%%%%%%%%%%%%%%%%%%%%%%%%%%%%%%%%%%%%%

\subsection{Quadrature Methods -- the Questions}

Given an integrand $ f : [ a, b ] \mapsto {\mathbb R} $, consider the
numerical approximation of the definite integral by a rule
$ \left\lbrace x_k, w_k \right\rbrace $ on $n$ points:
\begin{eqnarray*}
  \int_{a}^{b}
    f \left( x \right)
  dx
  \approx
  \sum_{k=1}^n
    f \left( x_k \right)
    w_k.
\end{eqnarray*}
The interval $ \left[ a, b \right] $ is possibly infinite or semi-infinite, but
this report considers only finite intervals; and without loss of
generality, let $ \left[ a, b \right] = \left[ 0, 1 \right] $. Similarly, the integrand
could include a weighting factor $ \omega \left( x \right) $, but this is not
required here.

The \emph{degree} of a quadrature rule is the maximal degree of the
polynomial that it can integrate exactly.\footnote{Comments on errors
refer to discretisation, not machine roundoff error unless explicitly
stated.} That is, if the degree of a rule on $n$ points is $p$,
then:
\begin{eqnarray*}
  \int_0^1
    x^j
  dx
  =
  \sum_{k=1}^n
    x_k^j w_k,
  \qquad
  j = 0:p.
\end{eqnarray*}

If $f$ is smooth, the rate of convergence for $n$ point Gau{\ss}ian
quadrature is $ {\cal O} \left( {\rho}^n \right) $ (for some $\rho<1$),
and for the composite Simpson's rule it is
$ {\cal O} \left( n^{-4} \right) $.  That is, the error decreases more
quickly for Gau{\ss}ian quadrature.  This is \emph{not} true in general
if $f$ has a singularity.%
\footnote{%
  `Singularity' is intended to always mean `end point singularity'. If
  a particular singularity in the integrand is not at an end point,
  then the interval can be subdivided so that the singularity is at the
  end points of the two subintervals.  Most quadrature methods perform
  poorly on integrands with internal singularities.
}
For example, consider the `square root' singularity
$ f ( x ) = \sqrt{x}, \; x \in \left[ 0, 1 \right] $. In this case, the
rate of convergence for Gau{\ss}ian quadrature falls to
$ {\cal O} \left( n^{-3} \right) $, whilst that of Simpson's rule is
$ {\cal O} \left( n^{-3/2} \right) $. Even so, using a \emph{composite}
Simpson's rule and a \emph{graded mesh}, a convergence rate of
$ {\cal O} \left( n^{-4} \right) $ can be recovered.

A \emph{composite} quadrature rule is created by subdividing the
interval of integration into $m$ subintervals, and evaluating the
integral over each subinterval using an appropriate quadrature rule.
Choosing $ x_{j-1} < x_j, \; j = 1:m $, and $ x_0 = 0, \; x_m = 1 $:
\begin{eqnarray*}
  \int_0^1
    f \left( x \right)
  dx
  =
  \sum_{j=1}^m
    \int_{x_{j-1}}^{x_j}
      f \left( x \right)
    dx.
\end{eqnarray*}

\vfill

\pagebreak

\emph{Grading} the mesh means that the subdivision is organised in some
way.  A description of a generalised mesh grading to cater for
complicated possibilities, such as \emph{adaptive} quadrature is found
in~\cite{ChandlerGraham:88}. Here, the simplification of nonadaptive
meshes is used. Meshes are graded by assigning mesh points according to
some simple function. A quadrature rule of degree $ p_j $ (a function
of $ n_j $, the number of points used by the rule) is used on the
interval $ \left[ x_{j-1}, x_j \right] $. This raises two issues:
\begin{itemize}
\item
  How should $ x_j $ and $ n_j $ be chosen? That is, how should
  the mesh be graded, and how should the degree of the
  quadrature rules on each subinterval vary?
\item
  How is this procedure dependent on the integrand? Consider the
  generalisation of $ \sqrt{x} $ to $ {| x |}^{\alpha} $, or
  more pathological cases. For experimental work, a known
  integrand is easy to deal with. What of more general cases,
  where no explicit functional information is known?
\end{itemize}
The remainder of this section describes some partial answers to these
questions, and displays some experiments. The answers deal with the
case $ {| x |}^{\alpha}, \; -1 < \alpha < 1 $, and the experiments
demonstrate the case $ \alpha = 1/2 $.

%%%%%%%%%%%%%%%%%%%%%%%%%%%%%%%%%%%%%%%%%%%%%%%%%%%%%%%%%%%%%%%%%%%%%%%%

\subsection{Gau{\ss}ian Quadrature}

Gau{\ss}ian quadrature rules are the \emph{best} possible rules in the
sense that they are of maximal degree.%
\footnote{%
  This may of course, not be ideal for the particular application, but
  in the absence of theoretical functional information about the
  integrand, nothing beyond this can be said about the convergence
  rates of \emph{any} quadrature method. In practice, with commonly
  occurring functions, there is a certain amount of implicit
  theoretical information which can be used in error analysis.
}
This is due to the exploitation of the maximal number
of degrees of freedom in the choice of their nodes and weights.
Gau{\ss} methods divide into categories depending on the associated
weight function, and whether there are any prescribed quadrature
points. For the applications in this report, it is preferable to use
the end points of the interval as quadrature points, and a unit weight
function is assumed. The appropriate set of rules are called the
Gau{\ss}--Lobatto rules~\cite[pages 101--104]{DavisRabinowitz:84}.

\begin{theorem}[Gau{\ss}--Lobatto Quadrature]
  ~\\
  Given $ f \in C^{2n-2} \left[ a, b \right] $, the $n$ point
  Gau{\ss}--Lobatto quadrature rule ($ n \geqslant 2 $), has nodes
  $ a \equiv x_1 < x_2 < \dots < x_{n-1} < x_n \equiv b $,
  and positive weights $ w_1, \dots, w_n $ such that:
  \begin{eqnarray*}
    \int_a^b
      f \left( x \right)
    dx
    =
    \sum_{k=1}^n
      f \left( x_k \right)
      w_k
    +
    E_n.
  \end{eqnarray*}
  Here $ E_n $ is dependent on $f$, $a$, $b$ and $n$:
  \begin{equation}
    E_n
    =
    -
    \displaystyle\frac
    {
      n \left( n + 1 \right)
      {\left[ \left( n - 2 \right)! \right]}^4
    }{
      \left( 2n - 1 \right)
      { \left[ \left( 2n - 2 \right)! \right]}^3
    }
    {\left( b - a \right)}^{2n - 1}
    f^{ \left( 2n - 2 \right)} \left( \xi \right),
    \qquad
    \xi \in \left( a, b \right).
    \label{eq:pony}
  \end{equation}
  The rule is of degree $ p = 2n - 3 $ (this is always odd).
  Observe that for $ n = 2 $ the rule is the trapezoidal rule,
  and for $ n = 3 $ it is Simpson's rule.
  Only for $ n \geqslant 4 $ do these rules diverge from the
  series of closed Newton--Cotes rules (see Table \ref{tab:degrees} on
  page \pageref{tab:degrees}).
\end{theorem}
An algorithm~\cite{Golub:73} for finding
$ \left\lbrace x_k, w_k \right\rbrace $ using a matrix eigenvalue
technique is implemented in Appendix \ref{app:lobatto}.

\vfill

\pagebreak

%%%%%%%%%%%%%%%%%%%%%%%%%%%%%%%%%%%%%%%%%%%%%%%%%%%%%%%%%%%%%%%%%%%%%%%%

\subsection{Graded Meshes}

A number of different methods for grading meshes appear in the
literature. Three important methods~\cite{PostellStephan:90} are
described below. In each case, the mesh subdivides the interval
$ \left[ 0 , 1 \right] $, when the integrand has a singularity at $ 0 $.  The
meshes have $m$ subintervals, that is $ m + 1 $ points, including the
ends. The $j$th mesh point is at $ x_j $, and the $j$th interval is
of width $ h_j $:
\begin{enumerate}
\item
  \emph{Quasiuniform}. The mesh is essentially uniform; that is
  for some constant $ \tau < 1 $,
  $ h_j \in \left[ h \tau, h \right], \; j = 1:m-1 $, where
  $ h = \max_j \left\lbrace h_j \right\rbrace $.
\item
  \emph{Algebraic}. For some $ \gamma \geqslant 1 $,
  $ x_j = {\left( \displaystyle\frac{j}{m} \right)}^{\gamma} $, $ j = 0:m $.
\item
  \emph{Geometric}. For some $ 0 < \sigma < 1 $,
  $ x_j = {\sigma}^{m - j} $, $ j = 1:m $; and $ x_0 = 0 $.
\end{enumerate}
%%% The following purges deal with the loss of the meshpts.ps file ...
%%% These are illustrated in Figure \ref{fig:meshpts}.

A geometrically graded mesh is
%%% also
illustrated on one segment of a closed contour in Figure
\ref{fig:nomen} on page \pageref{fig:nomen}.

%%% \begin{figure}[htbp]
%%%   \begin{center}
%%%     \input{graphics/meshpts.ps}
%%%     \caption{Various gradings of mesh points.}
%%%     \label{fig:meshpts}
%%%   \end{center}
%%% \end{figure}

%%%%%%%%%%%%%%%%%%%%%%%%%%%%%%%%%%%%%%%%%%%%%%%%%%%%%%%%%%%%%%%%%%%%%%%%

\subsection{$h$, $p$ and $h$-$p$ Quadrature Methods}

The `$h$-$p$' nomenclature presented here originated in
papers by Babu\v{s}ka et al.~\cite{%
  BabuskaDorr:81,%
  GuiBabuska:86a,%
  GuiBabuska:86b,%
  GuiBabuska:86c%
}, on finite element methods, based on previous work which did not
explicitly use this schema. The following discussion of the three
methods refers to their use with \emph{graded} meshes.

%%%%%%%%%%%%%%%%%%%%%%%%%%%%%%%%%%%%%%%%%%%%%%%%%%%%%%%%%%%%%%%%%%%%%%%%

\subsubsection{$h$ Methods}

An $h$ quadrature method is composed using two steps:
\begin{enumerate}
\item
  Choose an underlying mesh of subintervals; possibly a graded
  mesh determined by the user, from analysis of the singularities
  of the integrand.
\item
  Integrate over each of the $m$ mesh intervals, applying
  the same $n$ point quadrature rule. The result is a
  composite quadrature rule on a total of $N$ points. These
  $N$ points will be called the \emph{node} points from now on.
  The functional relationship $ N \left( m, n \right) $ is dependent on
  whether the basic quadrature rule is open or
  closed. (If the rule is open, the original mesh points are not
  included in the final rule.) Observe that the user cannot
  arbitrarily select $N$, only $m$ and $n$.
  \begin{eqnarray*}
    N \left( m, n \right)
    =
    \left\lbrace
    \begin{array}{ll}
      m \left( n - 2 \right) + m + 1 &
      (\mathrm{closed~basic~rule})
      \\
      m n &
      (\mathrm{open~basic~rule}).
    \end{array}
    \right.
  \end{eqnarray*}
\end{enumerate}
A particular basic rule is decided upon (e.g. Simpson's rule), and
desired accuracy is hopefully attained by simply increasing $m$ (that
is, $N$). Whatever grading is chosen, the separation of the node
points ($h$) decreases as $m$ is increased, hence the name `$h$
method'. For a uniform mesh (which works well for smooth integrands),
$ h = \left( b - a \right) / \left( N - 1 \right) $ is constant.  Alternatively, open
rules, or Gau{\ss} rules can be used, the only important factor is that
all the basic rules are of the same type and degree.

%%%%%%%%%%%%%%%%%%%%%%%%%%%%%%%%%%%%%%%%%%%%%%%%%%%%%%%%%%%%%%%%%%%%%%%%

\subsubsection{$p$ Methods}

In a $p$ method, again a graded mesh is created. Integration is
performed over each mesh interval using a basic quadrature rule on $n$
points. Here, $n$ instead of $m$ is varied by the user.  That is, for a
given number of mesh subintervals, $m$, a set of rules of increasing
degree (that is $n$, the number of points involved) is used, until
desired accuracy is obtained.  The same functional relationship
$N\left(m,n\right)$ exists.  As $n$ is increased, the rules used grow
in their degree ($p$), hence the name `$p$ method' (see also Table
\ref{tab:degrees}).

To illustrate, consider the family of closed Newton--Cotes
rules. Assume that the interval has been subdivided, possibly using an
adaptive algorithm that chooses smaller subdivisions where there the
integrand has greater derivative. Approximate the integral over each
division using the trapezoidal rule ($ p = 1 $), and inspect the
result. If it is unacceptable, repeat using Simpson's rule ($ p = 3 $).
Continue this process until results are acceptable.

%%%%%%%%%%%%%%%%%%%%%%%%%%%%%%%%%%%%%%%%%%%%%%%%%%%%%%%%%%%%%%%%%%%%%%%%

\subsubsection{$h$-$p$ Methods}
\label{sec:gmfq}

The $h$-$p$ method is the natural combination of the two previous
methods. The user may vary both $m$ and $n$.  The idea behind this
is to create a composite rule that minimises errors in the
approximation, for a given number of node points $N$. (Experiment
demonstrates that this is achievable.) The user chooses a family
of basic quadrature rules, then decides how to vary $n$ with mesh
interval. As the singularities considered are \emph{always} at end
points, a good choice is to organise small mesh intervals and low
degree rules (small $n$) near the end points, and larger mesh
intervals and high degree rules away from them, where the integrand is
expected to be smooth.

A simple choice is to begin with a rule on $ n = 2 $ points on the
smallest interval, then linearly increase $n$ with the number of the
mesh interval. Other discrete integer functions $ n_j, \; j = 1:m $ are
easily designed. The only constraint on these functions is that if
any error analysis is to be done, there should be some regularity in
$ n_j $. (Choosing basic rules from the same family facilitates this.)
This implementation uses the Gau{\ss}--Lobatto rule of degree $ 1 $
($ n = 2 $) on the first interval, degree $ 3 $ ($ n = 3 $) on the
second, etc.

Creating the composite quadrature rule is quite difficult. Each of the
basic quadrature rules must be appropriately scaled and shifted, and
then coincident mesh points must be combined.  This is further
complicated in the cases of closed meshes, closed quadrature rules, and
contour integration, where the end points of various segments of the
parameterisation must also be combined. (This is exacerbated if the
contour is closed.) The CBIEM requires all of these to be implemented.
The (closed) contours involved have corners, and the integrand will
usually have singularities at these corners. As it will be important to
keep the collocation points (see \S\ref{sec:cbiedisc}) between, and not
on, the corners, the underlying meshes \emph{must} include end points.
This means that the basic quadrature rules must be closed, so as to
include the end points.

The literature recommends using a geometrically graded mesh, with
an $h$-$p$ quadrature method. (This is implemented in the CBIEM.)
For maximum efficacy, the basic quadrature rules chosen must be of
maximal degree, which restricts them to Gau{\ss} rules. They must also
be closed. An $n$ point rule already has two of its points fixed, at
the ends. The appropriate rule is known as the Gau{\ss}--Lobatto rule,
which is of degree $ p = 2n - 3 $.

%%%%%%%%%%%%%%%%%%%%%%%%%%%%%%%%%%%%%%%%%%%%%%%%%%%%%%%%%%%%%%%%%%%%%%%%

\subsection{Error Analysis for the $h$ and $h$-$p$ Methods}

This section is tedious, and consists mainly of technical arguments.
The important parts are Theorem \ref{thm:horder} on page
\pageref{thm:horder}; Theorem \ref{thm:hporder} on page
\pageref{thm:hporder}; and the experimental results in 
\S\S\ref{sec:realint} and \ref{sec:ccontint}. The rest can be skipped
without loss of continuity.

%%%%%%%%%%%%%%%%%%%%%%%%%%%%%%%%%%%%%%%%%%%%%%%%%%%%%%%%%%%%%%%%%%%%%%%%

\subsubsection{Error Analysis for the $h$ Method}

This section computes an error bound for the $h$ method using an
algebraic mesh, for the integrand $ {| x |}^{\alpha} $, on
$\left[ 0, 1 \right] $.  Clearly $ \alpha > -1 $ is necessary for the
integral to be proper, and thus make its computation sensible.  For
$ -1 < \alpha < 0 $, $ {| x |}^{\alpha} \notin C^0 \left[ 0, 1 \right]$,
so the integrand is unbounded, but the integral is nonetheless defined.
If $ 0 \leqslant \alpha < 1 $,
${|x|}^{\alpha} \in C^0 \left[ 0, 1 \right] $, but
${|x|}^{\alpha} \notin C^1 \left[ 0, 1 \right] $.
If $ \alpha \in \mathbb{N} $, then $ {| x |}^{\alpha} \in C^{\infty}
\left[ 0, 1 \right] $, so the singularity vanishes and the case is of
lesser interest.  If $ \alpha > 1 $, but $ \alpha \notin \mathbb{N} $,
then all of the higher derivatives at $ x = 0 $ will not exist. Using
the notation $ \lfloor x \rfloor $ and $ \lceil x \rceil $ for the
least integers (respectively) greater than and less than $ x \in
{\mathbb R} $, in general, for $ \alpha > 1 $, $ {| x |}^{\alpha} \in
C^{\lfloor \alpha \rfloor} \left[ 0, 1 \right] $, but $ {| x
|}^{\alpha} \notin C^{\lceil \alpha \rceil} \left[ 0, 1 \right] $.
These cases are not particularly interesting, so the limit $ \alpha < 1
$ is made for simplicity. That is, consider $ -1 < \alpha < 1 $, which
includes the paradigm example $ x^{1/2} $. The $n$th derivative of $ f
\left( x \right) = {| x |}^{\alpha} $, for $ x \neq 0 $, is:

\begin{eqnarray*}
  f^{(n)} \left( x \right)
  =
  \displaystyle\frac
  {
    \Gamma \left( \alpha + 1 \right)
  }{
    \Gamma \left( \alpha + 1 - n \right)
  }
  {| x |}^{\alpha - n}.
\end{eqnarray*}

Consider the interval $ \left[ 0, 1 \right] $ partitioned into $m$
subintervals, where $ x_j $ is the $j$th mesh point, $ j = 0:m $, and
$ h_j = x_j - x_{j-1} $ is the width of the $j$th interval, for
$ j = 1:m $. Recall, for an \emph{algebraic} grading, a real constant
$ \gamma \geqslant 1 $ is chosen,%
\footnote{%
  Choosing $ \gamma < 1 $ results in some mesh points possibly lying
  outside the interval $ \left[ 0, 1 \right] $.
}
and the mesh is defined by:
$ x_j = {\left( j/m \right)}^{\gamma} $, $ j = 0:m $. Differentiating
and applying the mean value theorem shows:
\begin{eqnarray*}
  h_j
  =
  {\left( \displaystyle\frac{j}{m} \right)}^{\gamma}
  -
  {\left( \displaystyle\frac{j-1}{m} \right)}^{\gamma}
  =
  \displaystyle\frac{1}{m^{\gamma}}
  \left[
    j^{\gamma}
    -
    {\left( j - 1 \right)}^{\gamma}
  \right]
  \leqslant
  \displaystyle\frac{\gamma}{m}
  {\left( \displaystyle\frac{j}{m} \right)}^{\gamma-1}
  =
  {\displaystyle \frac{d x_j}{d j}}.
\end{eqnarray*}
For a \emph{geometric} grading, for some constant
$ 0 < \sigma < 1 $, $ x_j = {\sigma}^{m-j} $, $ j = 1:m $; $ x_0 = 0 $,
so:
\begin{eqnarray*}
  h_j
  =
  {\sigma}^{m-j}
  -
  {\sigma}^{m-j+1}
  =
  \left( 1 - \sigma \right)
  {\sigma}^{m-j},
  \qquad
  j = 2:m
  ,
  \qquad
  h_1 = {\sigma}^{m-1}.
\end{eqnarray*}
Consider an algebraic grading, using closed basic quadrature rules.
For an $h$ or $h$-$p$ method, the integral on the $j$th
interval is computed using a quadrature rule with
$ n_j \geqslant 2 $ points. For an $h$ method, $ n_j $ is constant;
for instance, $ n_j \equiv 2 $ means that the integral over each mesh
interval is computed using the trapezoidal rule -- the rule is a
\emph{composite} trapezoidal rule. The degrees of some common quadrature
rules are presented in Table \ref{tab:degrees}.

\begin{table}[htbp]
  \centering
  \begin{tabular}{||c||c|l||}
    \hline\hline
    Degree & Gau{\ss}--Lobatto & Newton--Cotes \\
    \hline\hline
    1  & 2  & 2 (trapezoidal) \\
    \hline
    3  & 3  & 3 (Simpson's), 4 (3/8ths) \\
    \hline
    5  & 4  & 5 (Boole's), 6 \\
    \hline
    7  & 5  & 7, 8 \\
    \hline\hline
  \end{tabular}
  \caption
  {
    The number of points $n$, and associated degrees of
    some closed quadrature rules. For Gau{\ss}--Lobatto
    rules, the degree is $ 2n - 3 $, and for
    Newton--Cotes rules, the degree is $n$ if $n$ is
    odd, else it is $ n - 1 $.
  }
  \label{tab:degrees}
\end{table}

Now consider the global error $ E_m $ of the difference between the
true solution $I$ and the approximation $ I_m $, induced by the
quadrature on $m$ subintervals, that is $ I = I_m + E_m $. For the
function $ {| x |}^{\alpha} $, bounds on $ E_m $ are readily
found. Define $ e_j $ as the component of $ E_m $ due to the $j$th
mesh interval, that is $ E_m = \sum_{j=1}^m e_j $.  To bound $ E_m $,
note $ E_m \leqslant \sum_{j=1}^m | e_j | $, and then bound each of the
$ | e_j | $.

\vfill

\pagebreak

The error result for a degree $p$ quadrature rule on an interval
$ \left[ x_{j-1}, x_j \right] $, of width $ h_j $, with a function
$ f \left( x \right) \in C^{p+1} \left[ x_{j-1}, x_j \right] $, is:
\begin{equation}
  e_j
  =
  C \left( p \right)
  h_j^{p+2}
  f^{\left( p+1 \right)} \left( \xi \right),
  \qquad
  \xi \in \left[ x_{j-1}, x_j \right].
  \label{eq:turkey}
\end{equation}
For $n$ point Gau{\ss}--Lobatto quadrature, the degree is
$ p = 2n - 3 $. (\ref{eq:turkey}) is derived from (\ref{eq:pony}) by a
scaling argument, and:
\begin{equation}
  C \left( p \right)
  =
  -
  \displaystyle\frac
  {
    \left( p + 3 \right) \left( p + 5 \right)
    {\left[ \left( \left( p - 1\right) / 2 \right)! \right]}^4
  }{
    2^2
    \left( p + 2 \right)
    { \left[ \left( p + 1 \right)! \right]}^3
  }.
  \label{eq:ram}
\end{equation}
As $ {| x |}^{\alpha} \in C^{\infty} \left( 0, 1 \right] $, (\ref{eq:turkey}) holds
for every mesh interval except the first, where the error is known
exactly, e.g. for the trapezoidal rule:
\begin{equation}
  e_1
  =
  \displaystyle\frac{h_1^{\alpha + 1}}{\alpha + 1}
  -
  \frac{h_1}{2} h_1^{\alpha}
  =
  \left(
    \displaystyle\frac{1}{\alpha + 1}
    -
    \displaystyle\frac{1}{2}
  \right)
  h_1^{\alpha + 1}.
  \label{eq:goat}
\end{equation}
As the maximum value of the $ \left( p + 1 \right) $th derivative of
$ {| x |}^{\alpha} $ on the $j$th interval, is at its left hand end,
$ f^{\left( p+1 \right)} \left( \xi \right) \leqslant f^{\left( p+1 \right)} \left( x_{j-1} \right) $,
the total error can be bounded:
\begin{eqnarray*}
  E_m
  \leqslant
  | e_1 |
  +
  C
  \sum_{j=2}^m
    h_j^{p+2}
    f^{\left( p+1 \right)} \left( x_{j-1} \right)
  \leqslant
  | e_1 |
  +
  C
  \sum_{j=1}^{m-1}
    h_{j+1}^{p+2}
    f^{\left( p+1 \right)} \left( x_j \right).
\end{eqnarray*}
Here, $C$ refers to a positive constant, independent of $m$, that
may vary from line to line.  Substituting for $ h_{j+1} $ and
$ f^{\left( p+1 \right)} \left( x_j \right) $ using an algebraically
graded mesh yields:
\begin{eqnarray*}
  E_m
  & \leqslant &
  | e_1 |
  +
  C
  \sum_{j=1}^{m-1}
    {\left[
      \displaystyle\frac{\gamma}{m}
      {\left(
        \displaystyle\frac{j-1}{m}
      \right)}^{\gamma-1}
    \right]}^{p+2}
    \displaystyle\frac
    {
      \Gamma \left( \alpha + 1 \right)
    }{
      \Gamma \left( \alpha - p \right)
    }
    {\left(
      \displaystyle\frac{j}{m}
    \right)}^{\gamma \left( \alpha - p - 1 \right)}
  \\
  & \leqslant &
  | e_1 |
  +
  C
  \sum_{j=1}^{m-1}
    {\left[
      \displaystyle\frac{\gamma}{m}
      {\left(
        \displaystyle\frac{j}{m}
      \right)}^{\gamma-1}
    \right]}^{p+2}
    {\left(
      \displaystyle\frac{j}{m}
    \right)}^{\gamma \left( \alpha - p - 1 \right)}.
\end{eqnarray*}
As $ | e_1 | $ is less than some constant multiplied by the `$m$th
term' in the sum:
\begin{equation}
  E_m
  \leqslant
  \displaystyle\frac
  {
    C
  }{
    m^{\gamma \left( \alpha + 1 \right)}
  }
  \sum_{j=1}^m
    j^{\gamma \left( \alpha + 1 \right) - \left( p+2 \right)}.
  \label{eq:parrot}
\end{equation}

\pagebreak

Simplification of (\ref{eq:parrot}) (see below) leads to the result:
\begin{theorem}[Convergence of the $h$ Method with Algebraic Grading]
  ~\\
  Consider the approximation of $ \int_0^1 {| x |}^{\alpha} dx $,
  $ -1 < \alpha < 1 $, using an $h$ method based on a
  quadrature rule of degree $p$ (an odd positive integer), on
  an algebraic mesh on a total of $ m \geqslant 2 $ intervals,
  with mesh parameter $ \gamma \geqslant 1 $. For some constant
  $C$, the error $ E_m $ satisfies:
  \begin{eqnarray*}
    E_m
    \leqslant
    C m^{-z}.
  \end{eqnarray*}
  Here $z$ is:
  \begin{eqnarray*}
    z
    =
    \left\lbrace
    \begin{array}{ll}
      \gamma \left( \alpha + 1 \right)
      &
      1 \leqslant \gamma <
      \left( p + 1 \right) / \left( \alpha + 1 \right)
      \\
      p + 1
      &
      \mathrm{else}.
    \end{array}
    \right.
  \end{eqnarray*}
  When $ \gamma = \left( p + 1 \right) / \left( \alpha + 1 \right) $,
  $ E_m \leqslant C \ln \left( m \right) / m^{p+1} $.
  \label{thm:horder}
\end{theorem}
\begin{proof}
  Take (\ref{eq:parrot}) and write $ E_m $ as:
  \begin{eqnarray*}
    E_m
    & \leqslant &
    \displaystyle\frac
    {
      C
    }{
      m^{\gamma \left( \alpha + 1 \right)}
    }
    \sum_{j=1}^m
      j^{\gamma \left( \alpha + 1 \right) - \left( p+2 \right)}
    \leqslant
    \displaystyle\frac
    {
      C
    }{
      m^{\gamma \left( \alpha + 1 \right)}
    }
    \int_1^m
      x^{\gamma \left( \alpha + 1 \right) - \left( p+2 \right)}
    dx
    \\
    & \leqslant &
    \displaystyle\frac
    {
      C
    }{
      m^{\gamma \left( \alpha + 1 \right)}
    }
    \int_1^{\infty}
      x^{\gamma \left( \alpha + 1 \right) - \left( p+2 \right)}
    dx.
  \end{eqnarray*}
  The integral converges if:
  $
    \gamma \left( \alpha + 1 \right) - \left( p+2 \right) < - 1
  $, and in this case, it converges to the constant
  $ {\left[ \left( p+1 \right) - \gamma \left( \alpha + 1 \right) \right]}^{-1} $,
  independent of $m$. Absorbing this into the main constant,
  the result for
  $ 1 \leqslant \gamma < \left( p+1 \right) / \left( \alpha+1 \right) $ is created.

  In the case $ \gamma = \left( p + 1 \right) / \left( \alpha + 1 \right) $,
  (\ref{eq:parrot}) is
  \begin{eqnarray*}
    E_m
    \leqslant
    \displaystyle\frac{C}{m^{p + 1}}
    \sum_{j=1}^m
      \displaystyle\frac{1}{j}
    \leqslant
    \displaystyle\frac{C}{m^{p + 1}}
    \int_1^m
      \displaystyle\frac{1}{x}
    dx
    \leqslant
    \displaystyle\frac{C}{m^{p + 1}}
    \ln \left( m \right).
  \end{eqnarray*}
  Thus again, the error is bounded by a constant, this time,
  dependent on $m$. Observe that
  $ p+1 $ is immediately able to be replaced with
  $ \gamma \left( \alpha + 1 \right) $, demonstrating the continuity of
  the formulae.

  Lastly, bound $ E_m $ as:
  \begin{eqnarray*}
    E_m
    & \leqslant &
    \displaystyle\frac{C}{m^{p + 2}}
    \sum_{j=1}^m
      {\left(
        \displaystyle\frac{j}{m}
      \right)}^{
        \gamma \left( \alpha + 1 \right) - \left( p+2 \right)
      }
    \leqslant
    \displaystyle\frac{C}{m^{p + 2}}
    \int_{1/m}^1
      x^{\gamma \left( \alpha + 1 \right) - \left( p+2 \right)}
    dx
    \\
    & \leqslant &
    \displaystyle\frac{C}{m^{p + 2}}
    \int_0^1
      x^{\gamma \left( \alpha + 1 \right) - \left( p+2 \right)}
    dx.
  \end{eqnarray*}
  As $ \gamma \left( \alpha + 1 \right) - \left( p + 2 \right) > - 1 $,
  the integral converges to
  $ {\left[ \gamma \left( \alpha + 1 \right) - \left( p + 1 \right) \right]}^{-1} $, again,
  another constant independent of $m$, which is absorbed into
  $C$. The second result is thus achieved, by observing
  that $ m^{-\left(p+2\right)} < m^{-\left(p+1\right)} $.
\end{proof}
The moral of this is that for a particular choice of $ \alpha $ and
$p$, there is an ideal choice of $ \gamma $, that is
$ {\gamma}^{*} = \left( p + 1 \right) / \left( \alpha + 1 \right) $, beyond which the
order of the error will not decrease. (Choosing $ \gamma $ greater than
this may reduce $C$.) Varying $ \gamma $ makes no difference to
computational expense.%
\footnote{%
  Choosing $ \lceil {\gamma}^{*} \rceil $ may be cheaper than
  nonintegral choices of $ {\gamma}^{*} $.
}
Observe that $ {\gamma}^{*} $ becomes unbounded as $\alpha\to -1^+$.
This is not surprising, as the integral itself becomes unbounded.  A
similar result to Theorem \ref{thm:horder} exists for a geometric
mesh.

\vfill

In summary, using an algebraic mesh of $m$ subintervals, and an $h$
method with a degree $p$ quadrature rule on each mesh interval, can
achieve, for $ {| x |}^{\alpha} $ integrands, an
$ {\cal O} \left( m^{-(p+1)} \right) $ convergence rate.  This is a
significant improvement on the equivalent result for a uniform mesh,
which is $ {\cal O} \left( m^{-\left(\alpha + 1\right)} \right) $.  The
exponential convergence rate for smooth integrands is not achieved, but
can be, using an $h$-$p$ method. The methods may be extended to
integrands of the form $ {| x |}^{\alpha} f \left( x \right) $, for
smooth functions $f$.

%%%%%%%%%%%%%%%%%%%%%%%%%%%%%%%%%%%%%%%%%%%%%%%%%%%%%%%%%%%%%%%%%%%%%%%%

\subsubsection{Error Analysis for the $h$-$p$ Method}

This section discusses the expected order of the error for the
$h$-$p$ method, for algebraic or geometric meshes using integrands
with end point singularities. Again, consider the integral
$ \int_0^1 {| x |}^{\alpha} dx $. For the $h$ method, $ p_j $, the
degree of the quadrature rule on the $j$th interval, was constant.
For the $h$-$p$ method, it is a function of $j$.  A low degree
rule is used on the interval adjacent to the end point singularity; and
higher degree rules are used on intervals away from it.  A simple
choice is to use a rule on $ 2 $ points on the first interval, and
increase the number of points linearly with $j$. That is, using
Gau{\ss}--Lobatto rules, where $ p_j = 2 n_j - 3 $; choosing
$ n_j = j + 1 $ gives $ p_j = 2 j - 1 $.

Recall the error result from (\ref{eq:turkey}), where
$C_j = C \left( p_j \right)$ is given by (\ref{eq:ram}):
\begin{eqnarray*}
  \int_{x_{j-1}}^{x_j}
    f \left( x \right)
  dx
  -
  I_{p_j}
  =
  e_j
  =
  C_j
  h_j^{p+2}
  f^{\left( p+1 \right)} \left( \xi \right),
  \qquad
  \xi \in \left[ x_{j-1}, x_j \right].
\end{eqnarray*}
Use Stirling's formula for large $x$ to approximate the factorials in
$ C_j $, as $ \left( x - 1 \right)! = \Gamma \left( x \right) $:
\begin{eqnarray*}
  \Gamma \left( x \right)
  \sim
  \sqrt
  {
    \displaystyle\frac{2 \pi}{x}
  }
  {\left(
    \displaystyle\frac{x}{e}
  \right)}^{x}
  \left\lbrace
    1
    +
    \displaystyle\frac{1}{12 x}
    +
    \dots
  \right\rbrace.
\end{eqnarray*}
Applying the first term of this to (\ref{eq:ram}) gives an asymptotic
bound for large $j$:
\begin{eqnarray*}
  C_j
  & = &
  -
  \displaystyle\frac{ \left( j + 1 \right) \left( j + 2 \right)}{2j + 1}
  \displaystyle\frac
  {
    {\left[ \left( j - 1 \right)! \right]}^4
  }{
    {\left[ \left( 2j \right)! \right]}^3
  }
  \; \sim \;
  -
  \displaystyle\frac{2^7}{e^3}
  \sqrt{\pi j}
  {\left(
    \displaystyle\frac{e}{8j}
  \right)}^{2j+3}.
\end{eqnarray*}
Thus:
\begin{equation}
  | C_j |
  \; \leqslant \;
  C
  \displaystyle\frac{{\mu}^j}{j^{5/2 + 2j}},
  \qquad
  \mu = e^2 / 2^6.
  \label{eq:rat}
\end{equation}
$ E_m $ cannot be bounded directly using this approximation for
$ C_j $, as it is not a proper bound.  Heuristically,%
\footnote{%
  This is brought out by experiment.
}
for either a geometric or an algebraic mesh, the rapid convergence of
$C_j $ to $ 0 $ means that $ E_m $ is expected to be dominated by
$e_1$.  Consider a mesh on $m$ intervals, and an associated composite
quadrature rule on a total of $ N = m \left( m + 1 \right) / 2 + 1 $
points. For a geometric mesh $ h_1 = {\sigma}^{m-1} $, and for an
algebraic mesh $ h_1 = m^{-\gamma} $.  Using (\ref{eq:goat}), this
gives:

\begin{equation}
  | e_1 |
  \; \leqslant \;
  C
  h_1^{\alpha+1}
  \; \leqslant \;
  C
  \left\lbrace
  \begin{array}{ll}
    m^{-\gamma \left( \alpha+1 \right)} & \mathrm{algebraic}
    \\
    {\sigma}^{\left( m - 1\right) \left( \alpha+1 \right)} &
      \mathrm{geometric}.
  \end{array}
  \right.
  \label{eq:wallaby}
\end{equation}

\vfill

\pagebreak

The error for an algebraic mesh is polynomial, whilst the error for a
geometric mesh is exponential.  If the errors with increasing $ N
\approx m^2/2 $ are plotted for a geometric mesh, an error of the form
$
  E_m
  \leqslant
  C
  {\rho}^{\sqrt{N}}
$
is observed (see Figure \ref{fig:hpintegp}), for some $ \rho < 1 $.

These rough results prompt more rigorous examination of $ E_m $. Recall:
\begin{eqnarray*}
  | e_j |
  & \leqslant &
  C_j
  h_j^{2j+3}
  f^{\left( 2j+2 \right)} \left( x_{j-1} \right).
\end{eqnarray*}
For a geometric mesh, $ x_j = {\sigma}^{m-j} $, so
$ h_j < \left( 1 - \sigma \right) {\sigma}^{m-j} $. Using
$ f \left( x \right) = {| x |}^{\alpha} $, gives:
\begin{eqnarray*}
  f^{(n)} \left( x \right)
  =
  \displaystyle\frac
  {
    \Gamma \left( \alpha + 1 \right)
  }{
    \Gamma \left( \alpha + 1 - n \right)
  }
  {| x |}^{\alpha - n}.
\end{eqnarray*}
The Stirling asymptotic approximation for $ C_j $ can be converted
to a genuine bound by observing for $ x \in \mathbb{N} $:
\begin{eqnarray*}
  \sqrt{\displaystyle\frac{2 \pi}{x}}
  {\left( \displaystyle\frac{x}{e} \right)}^x
  & \leqslant &
  \Gamma \left( x \right)
  \; \leqslant \;
  C
  \sqrt{\displaystyle\frac{2 \pi}{x}}
  {\left( \displaystyle\frac{x}{e} \right)}^x.
\end{eqnarray*}
Recall from (\ref{eq:rat}), that with $ \mu = e^2 / 2^6 $, for some $C$
independent of $x$:
\begin{eqnarray*}
  | C_j |
  & \leqslant &
  C
  \displaystyle\frac{{\mu}^j}{j^{2j+5/2}}.
\end{eqnarray*}
A bound on the error for the $j$th interval is now:\footnote{Strictly
speaking, this only applies for $ j = 2:m $, as the result for $ C_j $
only holds for $ j = 2:m $. Application of the result in
(\ref{eq:wallaby}) allows $ | e_1 | $ to be bounded by a constant multiple
of this $ C_1 $.}
\begin{eqnarray*}
  | e_j |
  & < &
  C
  \displaystyle\frac
  {
    {\mu}^j
  }{
    j^{2j+5/2}
  }
  {\left[
    \left( 1 - \sigma \right)
    {\sigma}^{m-j}
  \right]}^{2j+3}
  f^{\left( 2j+2 \right)} \left( {\sigma}^{m-j-1} \right).
\end{eqnarray*}
Simplification of this leads to a bound on $ E_m $. Observe that
$ {\left( 1 - \sigma \right)}^{2j+3} < 1 $, so:
\begin{eqnarray*}
  | e_j |
  & < &
  C
  \displaystyle\frac
  {
    {\mu}^j
  }{
    j^{2j+5/2}
  }
  {\sigma}^{\left( m-j \right) \left( 2j+3 \right)}
  \displaystyle\frac
  {
    \Gamma \left( \alpha + 1 \right)
  }{
    \Gamma \left( \alpha - 2j - 1 \right)
  }
  {\sigma}^{\left( m-j-1 \right) \left( \alpha - 2j - 2 \right)}.
\end{eqnarray*}
Absorb $ \Gamma \left( \alpha + 1 \right) $ into $C$, and use the
Stirling approximation for $ \Gamma \left( \alpha - 2j - 1 \right) $:
\begin{eqnarray*}
  \displaystyle\frac{1}{\Gamma \left( \alpha - 2j - 1 \right)}
  & \leqslant &
  \sqrt{\displaystyle\frac{\alpha - 2j - 1}{2 \pi}}
  {\left(
    \displaystyle\frac{e}{\alpha-2j-1}
  \right)}^{\alpha-2j-1}.
\end{eqnarray*}
Rearranging the exponent of $ \sigma $, and absorbing the term
$ e^{\alpha-1} $ into $C$:
\begin{eqnarray*}
  | e_j |
  & < &
  C
  \displaystyle\frac{{\mu}^j}{j^{2j+5/2}}
  {\sigma}^{\left( m-j \right) \left( \alpha + 1 \right) - \alpha + 2j + 2}
  \sqrt{\displaystyle\frac{\alpha - 2j - 1}{2 \pi}}
  {\left(
    \displaystyle\frac{e}{\alpha-2j-1}
  \right)}^{\alpha-2j-1}
  \\
  & < &
  C
  \displaystyle\frac
  {
    {\mu}^j
    {\sigma}^{m \left( \alpha + 1 \right)}
    {\sigma}^{j \left(  1 - \alpha \right)}
    e^{-2j}
  }{
    j^{2j+5/2}
    {\left( \alpha-2j-1 \right)}^{\alpha-2j-1/2}
  }.
\end{eqnarray*}
Expand $ \mu = e^2 / 2^6 $, and observe that
$ {\left( \alpha-2j-1 \right)}^{\alpha-1/2} > 1 $:
\begin{eqnarray*}
  | e_j |
  & < &
  C
  \displaystyle\frac
  {
    e^{2j}
    {\sigma}^{m \left( \alpha + 1 \right)}
    {\sigma}^{j \left(  1 - \alpha \right)}
    e^{-2j}
  }{
    j^{2j+5/2} 2^{6j}
    {\left( \alpha-2j-1 \right)}^{-2j}
  }
  \; < \;
  C
  \displaystyle\frac
  {
    {\sigma}^{m \left( \alpha + 1 \right)}
    {\sigma}^{j \left(  1 - \alpha \right)}
  }{
    2^{6j}
    j^{5/2}
  }
  {\left(
    \displaystyle\frac
    {
      \alpha-2j-1
    }{
      j
    }
  \right)}^{2j}.
\end{eqnarray*}
As $ | \left( \alpha - 2j - 1 \right) / j | < 2 $, then:
\begin{eqnarray*}
  | e_j |
  & < &
  C
  \displaystyle\frac
  {
    2^{2j}
    {\sigma}^{m \left( \alpha + 1 \right)}
    {\sigma}^{j \left(  1 - \alpha \right)}
  }{
    2^{6j}
    j^{5/2}
  }
  \; < \;
  C
  \displaystyle\frac
  {
    {\sigma}^{m \left( \alpha + 1 \right)}
    {\sigma}^{j \left(  1 - \alpha \right)}
  }{
    2^{4j}
    j^{5/2}
  }.
\end{eqnarray*}
Combine these, to bound:
\begin{eqnarray*}
  E_m
  & \leqslant &
  \sum_{j=1}^m
    | e_j |
  \; < \;
  C
  \sum_{j=1}^m
    \displaystyle\frac
    {
      {\sigma}^{m \left( \alpha + 1 \right)}
      {\sigma}^{j \left(  1 - \alpha \right)}
    }{
      2^{4j}
      j^{5/2}
    }
  \; < \;
  C
  {\sigma}^{m \left( \alpha + 1 \right)}
  \sum_{j=1}^m
    \displaystyle\frac
    {
      {\sigma}^{j \left(  1 - \alpha \right)}
    }{
      2^{4j}
      j^{5/2}
    }.
\end{eqnarray*}
Observing that $ 2^{4j} j^{5/2} > 1 $, this can be written as:%
\footnote{This throws away a lot of information!}
\begin{eqnarray*}
  E_m
  \; < \;
  C
  {\sigma}^{m \left( \alpha + 1 \right)}
  \sum_{j=1}^m
    {\left[ {\sigma}^{\left( 1 - \alpha \right)} \right]}^j.
\end{eqnarray*}
Using the result for the sum of a geometric progression:
\begin{eqnarray*}
  \sum_{j=1}^m
    {\left[ {\sigma}^{\left( 1 - \alpha \right)} \right]}^j
  & = &
  \displaystyle\frac
  {
    {\left( {\sigma}^{1 - \alpha} \right)}^{m+1} - 1
  }{
    {\sigma}^{1 - \alpha} - 1
  }
  \; = \;
  \displaystyle\frac
  {
    {\sigma}^{\left( 1 - \alpha \right) \left( m+1 \right)} - 1
  }{
    {\sigma}^{1 - \alpha} - 1
  }
  \; < \;
  C
  {\sigma}^{\left( 1 - \alpha \right) \left( m+1 \right)}.
\end{eqnarray*}
Thus:
$
  E_m
  \; < \;
  C
  {\sigma}^{m \left( \alpha + 1 \right)}
    {\sigma}^{\left( 1 - \alpha \right) \left( m+1 \right)}
  \; < \;
  C
  {\sigma}^{2m}
$.
As $ N \approx m^2 / 2 $, using $ \rho = \sigma^{2\sqrt{2}} $, this
simplifies to:
\begin{eqnarray*}
  E_m
  \; < \;
  C
  {\rho}^{\sqrt{N}}.
\end{eqnarray*}
This constant $C$ is actually a function of $ \alpha $ and
$ \sigma $, and the fact that the degrees of the quadrature rules are
linearly graded. This proves Theorem \ref{thm:hporder}:

\begin{theorem}[Convergence of the $h$-$p$ Method
                with Geometric Grading]
  ~\\
  Consider the approximation of $ \int_0^1 {| x |}^{\alpha} dx $,
  $ -1 < \alpha < 1 $, using an $h$-$p$ method based on a
  geometric mesh on $m$ intervals with parameter $ \sigma $,
  and using a Gau{\ss}--Lobatto quadrature rule on $ j + 1 $
  points on interval $j$.  For some constants $C$, and
  $ \rho < 1 $, the error $ E_m $ satisfies:
  \begin{eqnarray*}
    E_m
    \; < \;
    C
    {\sigma}^{2m}.
  \end{eqnarray*}
  \label{thm:hporder}
\end{theorem}
Babu\v{s}ka et al.~\cite{%
  BabuskaDorr:81,%
  GuiBabuska:86a,%
  GuiBabuska:86b,%
  GuiBabuska:86c%
} describe an approximation theory for the $h$, $p$ and $h$-$p$
methods for the finite element method.  This material may be able to be
simplified and adapted (as the theory for integration should be easier
than for approximation), and also used in error analysis.

%%%%%%%%%%%%%%%%%%%%%%%%%%%%%%%%%%%%%%%%%%%%%%%%%%%%%%%%%%%%%%%%%%%%%%%%

\subsubsection{Experimental Results for a Real Integral}
\label{sec:realint}

\begin{figure}[htbp]
  \begin{center}
    \epsfig{file=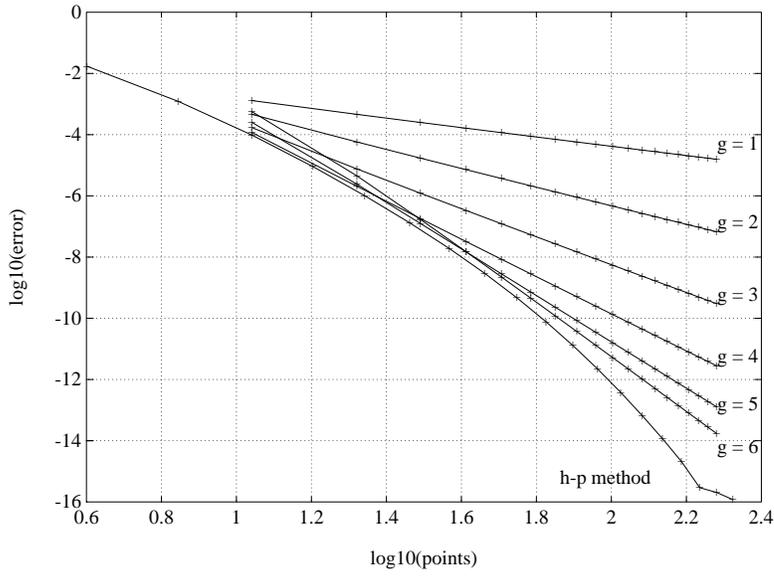,height=230pt}
    \caption{
      Errors for $h$ and $h$-$p$ methods applied to
      $ \int_0^1 \protect\sqrt{x} dx $, for various choices of
      $ \gamma \equiv g $. The slopes of the lines are approximately
      $ - 3 \gamma / 2 $. As the error tends to machine precision
      ($ \epsilon \approx 10^{-16} $), the convergence results lose
      their regularity.
    }
    \label{fig:hpintegg}
  \end{center}
\end{figure}

\begin{figure}[htbp]
  \begin{center}
    \epsfig{file=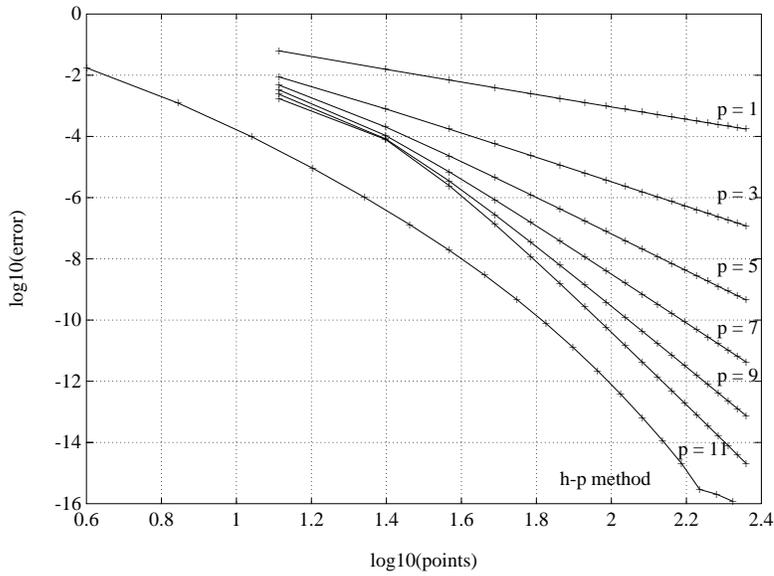,height=230pt}
    \caption{
      Errors for $h$ and $h$-$p$ methods applied to
      $ \int_0^1 \protect\sqrt{x} dx $, for constant $p$.
      The slopes of the lines are approximately
      $ - \left( p + 1 \right) $.
    }
    \label{fig:hpintegp}
  \end{center}
\end{figure}

The example $ \int_0^1 \sqrt{x} dx $ is used to demonstrate the above
convergence results. \textsc{matlab} code used (\texttt{hpmeth.m} and
\texttt{funchp.m}), is contained in Appendix \ref{app:code}, and error
results are presented in Figures \ref{fig:hpintegg} and
\ref{fig:hpintegp}.

Figure \ref{fig:hpintegg} compares convergence rates for various
choices of the algebraic grading parameter $ \gamma $, whilst holding
constant the number of points in the quadrature rule used on each
interval; that is $ p = 6 $ is fixed. Figure \ref{fig:hpintegp} shows
similar data, but varies $p$. (Here, $ \gamma $ is allowed to vary,
and is chosen to be equal to $p$ for convenience.) Both plots are
shown compared with the corresponding $h$-$p$ result, using a
geometric grading, and $ \sigma = 0.15 $. The linear increase in number
of points used in the quadrature rule means that $ n_j = j+1 $,
$ j = 1:m $, so $ N = m \left( m + 1 \right) / 2 + 1 $ (as rules are
closed, but ends are not), and thus $ m \propto \sqrt{N} $.  The
$h$-$p$ result demonstrates $ {\cal O} \left( {\rho}^{\sqrt{N}}\right)$
behaviour. The $h$-$p$ data does not show as a straight line, but this
can be seen in Figure \ref{fig:cint}, where it is plotted versus
$\sqrt{N}$.

\vfill

\clearpage

%%%%%%%%%%%%%%%%%%%%%%%%%%%%%%%%%%%%%%%%%%%%%%%%%%%%%%%%%%%%%%%%%%%%%%%%

\subsubsection{Extension to Complex Contour Integrals}
\label{sec:ccontint}

\begin{figure}[htbp]
  \begin{center}
    \epsfig{file=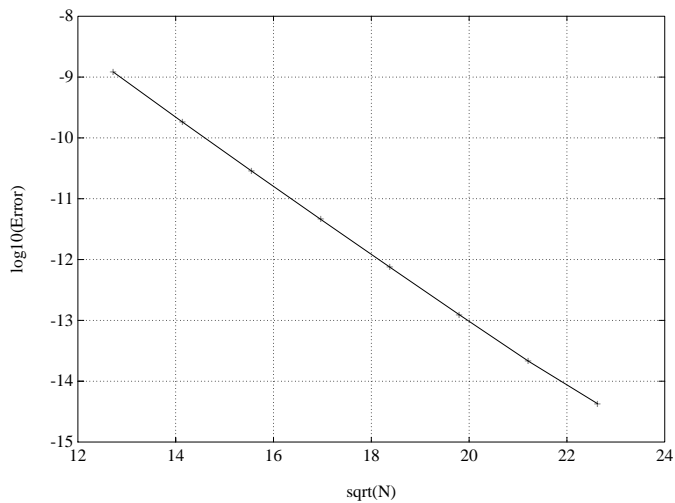,height=200pt}
    \caption{Error results for the complex contour integral.}
    \label{fig:cint}
  \end{center}
\end{figure}

The method is easily extended to complex contour integrals. Good test
problems have closed contours, and integrals which can be directly
evaluated using Cauchy's integral formula.  To demonstrate this,
consider $ f \left( z \right) = e^{-i \pi/4} z^{-1} \sqrt{z-1} $,
integrated around the unit circle. The integrand has a simple pole at
$z=0$, and a derivative singularity at $ z = 1 $. The resulting
integral is:
\begin{eqnarray*}
  \oint_{\Gamma}
    f \left( z \right)
  dz
  =
  e^{i \pi/4}.
\end{eqnarray*}
\textsc{matlab} code used (\texttt{cint.m} and \texttt{funcci.m}) is
contained in Appendix \ref{app:code}. Error results for an $h$-$p$
method using Gau{\ss}--Lobatto quadrature rules are presented in Figure
\ref{fig:cint}. The mesh is geometrically graded, with parameter
$\sigma = 0.15$. For a segment of a closed contour, with a corner at
either end, let $D$ be chosen as the number of mesh intervals between
each corner and a wide, central interval, so $ m = 2D + 1 $ is the
number of mesh intervals over that segment.  Here, as the contour is
the unit circle, $ 2 $ artificial corners are placed, and $D$ is varied
from $ 8 $ to $ 15 $.  The grading of the quadrature rules is similar
to that used in \S\ref{sec:realint} -- the number of points used in the
quadrature rule increases linearly with the number of mesh intervals
from the nearest corner, starting at $ 2 $ on the interval nearest the
corner, and finishing at $ D + 2 $ on the central interval.

Convergence is plotted for the logarithm of the error with $\sqrt{N}$.
Observe that the plot is linear, that is, the error is
$ {\cal O} \left( {\rho}^{\sqrt{N}} \right) $. These superb results
show that the method is excellent for the numerical approximation of
closed complex contour integrals.  This success motivates the use of
the $h$-$p$ method in the CBIEM, where quadrature rules for complex
contour integrals are required in the numerical approximation of the
solution to an integral equation.

%%%%%%%%%%%%%%%%%%%%%%%%%%%%%%%%%%%%%%%%%%%%%%%%%%%%%%%%%%%%%%%%%%%%%%%%

\subsection{Summary -- Advantages of $h$-$p$ Methods}

This section has discussed three important aspects of the numerical
approximation of integrals with end point singularities:

\begin{enumerate}
\item
  An $h$-$p$ quadrature method is superior to other methods.
\item
  A geometrically graded mesh is superior to other choices of
  grading (maybe only marginally better than an algebraic one).
\item
  The appropriate family of quadrature rules to use is the
  Gau{\ss}--Lobatto, as they are closed, and of maximal
  degree for the number of quadrature points used.
\end{enumerate}

The quadrature rule used in \S\ref{sec:CBIEM} is chosen in this manner.

\pagebreak

%%%%%%%%%%%%%%%%%%%%%%%%%%%%%%%%%%%%%%%%%%%%%%%%%%%%%%%%%%%%%%%%%%%%%%%%

\section{The Complex Boundary Integral Equation Method}
\label{sec:CBIEM}

\subsection{Origins and Description}

The CBIEM is a technique which numerically approximates the solution of
the Dirichlet problem.%
\footnote{%
  The space containing the functions approximating the solution of the
  Dirichlet problem is a Sobolev space, which is a generalisation of
  the Banach space of continuous functions to include functions with
  `weak derivatives'.
}
It reformulates the solution of the Dirichlet problem as the real part
of a function which can be found as the solution of a complex boundary
integral equation.  The solution of a discretised version of this
integral equation is then found using a collocation technique. Finally,
a discretisation of Cauchy's integral formula is used to approximate
the solution to the original problem at interior points, based on the
approximate boundary data.

It is related to the `Complex Variable Boundary Element
Method'~\cite{Hromadka:84}, which is a Galerkin version of the same
technique, using `hat' functions as a basis. (The collocation method
creates an approximation to the boundary data by interpolating from
known data, whilst the Galerkin constructs an approximation in terms of
a series of basis functions defined on segments of the boundary.) As
originally stated, the CVBEM only works on polygonal
domains,\footnote{The CVBEM has also been generalised to doubly
connected domains~\cite{KassabHsieh:90}.} whilst the CBIEM is more
general in that it also works on non-polygonal domains.

The Dirichlet problem considered is on an open, finite, simply
connected and non-empty region $ \Omega \subset \mathbb{R}^{2} $.
$\Omega$ is bounded by $ \Gamma $, a piecewise continuous,
anticlockwise oriented contour. $ \Gamma $ has a finite number of
corners, at which its derivative is discontinuous.  The Dirichlet
problem is:

Given boundary data $f$, find $ U : \Omega \mapsto \mathbb{R} $
subject to the conditions:
\begin{eqnarray*}
  {\nabla}^2 U \left( \mathbf{x} \right)
  =
  0,
  \qquad
  \mathbf{x} \in \Omega
  ,
  \qquad
  \qquad
  U \left( \mathbf{x} \right)
  =
  f \left( \mathbf{x} \right),
  \qquad
  \mathbf{x} \in \Gamma.
\end{eqnarray*}

Thus, the problem is to find the solution to Laplace's equation over
a region, given functional data around its perimeter. This has many
applications in the solution of potential problems, such as
electrostatics and fluid flow. The value of $U$ at points interior to
$ \Omega $ is determined by the boundary data being `diffused' from the
boundary inwards, according to the Laplacian operator.  It turns out
that the problem has a unique solution for all cases of $f$.  In all
but the most trivial cases, this solution is not expressible in closed
form, and a numerical approximation is required. With sufficient
(possibly enormous) computational effort, an approximation to any
degree of accuracy can usually be obtained.

Problems with `corner singularities' are of particular interest. In
these problems, $U$ is differentiable in the interior, but $ \nabla U
$ becomes unbounded as the corner is approached. It is known that this
behaviour is typical of solutions to the Dirichlet problem on domains
with corners. Even if the boundary data is smooth, $ \nabla U $ still
becomes singular near the corner. Numerical methods must be able to
produce good approximations to $U$, in spite of the corner
singularities.

`Interior' methods, such as finite difference and finite element
methods, become computationally expensive when applied to problems with
corner singularities, and boundary integral methods are more
appropriate. Interior methods require a finely discretised two
dimensional mesh in the region of the corner, which greatly increases
the size of the associated linear system. In contrast, a boundary
integral method has only to discretise its mesh in one dimension, that
of arc length on the boundary, and is expected to be much cheaper.

\vfill

\pagebreak

The usual boundary integral methods based on Green's functions lead to
a kernel with a logarithmic singularity, even on a smooth domain. This
is tedious to program, and computationally inefficient if high order
methods are used. If the CBIEM is used with singularity subtraction,
the integrands are smooth and can be done simply and accurately by
direct quadrature.

The problems caused by the corners and corner singularities are dealt
with using $h$-$p$ quadrature methods, and would be difficult to
implement with other types of integral equations~\cite{Chandler:90}.

\pagebreak

%%%%%%%%%%%%%%%%%%%%%%%%%%%%%%%%%%%%%%%%%%%%%%%%%%%%%%%%%%%%%%%%%%%%%%%%

\subsection{Development of the CBIEM}

The solution to the Dirichlet problem, $U$, is harmonic, as it
satisfies Laplace's equation in $ \Omega $. Identify
$ \mathbf{x} \in \mathbb{R}^{2} $ with $ z \in \mathbb{C} $. Now $U$ can
be thought of as the real component of an analytic function
$ W \left( z \right) = U \left( z \right) + i V \left( z \right) $,
where $V$ is uniquely determined to within a constant. $V$ can be
made unique by requiring $ V \left( {\zeta}_0 \right) = 0 $ for some
$ {\zeta}_0 \in \Gamma $ (see \S\ref{sec:solcoll}).  For all
$ z \in \Gamma $, $ U \left( z \right) \equiv f \left( z \right) $ is
immediately known.  The CBIEM first approximates $ V \left( z \right) $
on $ \Gamma $, and then uses Cauchy's integral formula to approximate
$W$, and hence $U$, at points within $ \Omega $.

%%%%%%%%%%%%%%%%%%%%%%%%%%%%%%%%%%%%%%%%%%%%%%%%%%%%%%%%%%%%%%%%%%%%%%%%

\subsubsection{Cauchy's Integral Formula}

For an analytic function $W$ on a bounded domain $ \Omega $,
Cauchy's integral formula is:
\begin{equation}
  \oint_{\Gamma}
    \displaystyle\frac
    {
      W \left( \zeta \right)
    }{
      \zeta - z
    }
  d\zeta
  =
  \pi i
  W \left( z \right)
  \times
  \left\lbrace
  \begin{array}{lll}
    0 & z \notin \Omega \cup \Gamma & \mathrm{exterior} \\
    1 & z \in \Gamma & \mathrm{boundary} \\
    2 & z \in \Omega & \mathrm{interior}.
  \end{array}
  \right.
  \label{eq:dolphin}
\end{equation}
When $z$ is on the boundary,\footnote{This result is a
simplification. If $z$ is at a corner, replace $ 1 $ with
$ \alpha / \pi $, where $ \alpha $ is the interior angle subtended by
the corner (else $ \alpha = \pi $). This result requires that
collocation points are \emph{not} placed at corners, to avoid unwanted
complexities in the implementation. Fortunately, this is already
overcome by the use of node points at the corners
(see \S\ref{sec:cbiedisc}).} the integral is a Hilbert transform. The
kernel is singular, and the result must be interpreted as a Cauchy
principal value integral~\cite[page 39]{CarrierKrookPearson:66}.  The
CBIEM requires approximation of the Cauchy integrals by quadrature.  In
\S\ref{sec:cbiedisc}, this is used to set up a linear system for the
approximation of $V$ on $ \Gamma $.  After this has been done, in
\S\ref{sec:singsub} it is used to compute an approximation to $W$
(and hence $U$) in the interior of $ \Omega $.

%%%%%%%%%%%%%%%%%%%%%%%%%%%%%%%%%%%%%%%%%%%%%%%%%%%%%%%%%%%%%%%%%%%%%%%%

\subsubsection{The Complex Boundary Integral Equation}
\label{sec:cbie}

To derive the integral equation underlying the CBIEM, observe that
letting $ W \left( \zeta \right) \equiv 1 $ for the case $z\in\Gamma$ in
(\ref{eq:dolphin}) yields:
\begin{eqnarray*}
  \oint_{\Gamma}
    \displaystyle\frac
    {
      1
    }{
      \zeta - z
    }
  d\zeta
  =
  \pi i,
  \qquad
  z \in \Gamma.
\end{eqnarray*}
Multiplying this by the constant $ W \left( z \right) $ gives:
\begin{eqnarray*}
  W \left( z \right)
  \oint_{\Gamma}
    \displaystyle\frac
    {
      1
    }{
      \zeta - z
    }
  d\zeta
  =
  \oint_{\Gamma}
    \displaystyle\frac
    {
      W \left( z \right)
    }{
      \zeta - z
    }
  d\zeta
  =
  \pi i W \left( z \right).
\end{eqnarray*}
Equating the $ \pi i W \left( z \right) $ with that in
(\ref{eq:dolphin}) gives:
\begin{equation}
  \oint_{\Gamma}
    \displaystyle\frac
    {
      W \left( \zeta \right) - W \left( z \right)
    }{
      \zeta - z
    }
  d\zeta
  =
  0,
  \qquad
  z \in \Gamma.
  \label{eq:seal}
\end{equation}

(\ref{eq:seal}) is called the `Complex Boundary Integral Equation'.
This derivation is parallel to that involved in
\emph{singularity subtraction}~\cite[page 184]{DavisRabinowitz:84}
and~\cite{Chandler:91}. The integrand is analytic, and as $\zeta\to z$,
it converges to $ W' \left( z \right) $. An analytic function $W=U+iV$,
which has as its real component the solution to the Dirichlet problem,
will satisfy (\ref{eq:seal}). The converse is also true -- a function
$W$ that satisfies (\ref{eq:seal}) will have a real component $U$ that
satisfies the Dirichlet problem. It is hoped that a function that
satisfies a discretisation of (\ref{eq:seal}) will have as its real
part the solution to a discretisation of the Dirichlet problem.

\vfill

\pagebreak

%%%%%%%%%%%%%%%%%%%%%%%%%%%%%%%%%%%%%%%%%%%%%%%%%%%%%%%%%%%%%%%%%%%%%%%%

\subsubsection{Discretisation of the CBIE}
\label{sec:cbiedisc}

\begin{figure}[htbp]
  \begin{center}
    \epsfig{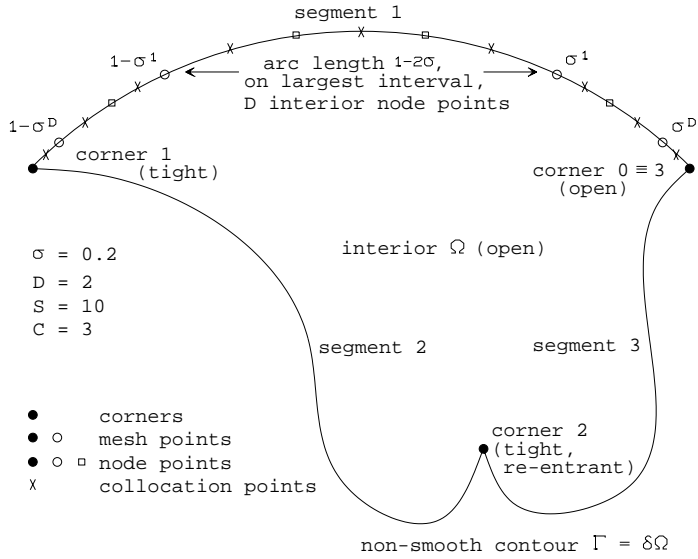}
    \caption{
      CBIEM nomenclature. Corner, mesh, node and collocation points on
      contour $ \Gamma $, about a region $ \Omega $. The lengths marked
      on segment $ 1 $ are the positions of geometric mesh points, in
      terms of a unit arc length on that segment. The node points
      correspond to a closed Newton--Cotes rule on each mesh interval.
    }
    \label{fig:nomen}
  \end{center}
\end{figure}

The CBIEM requires the numerical approximation of the (Cauchy) integral
in (\ref{eq:seal}), and this is achieved using quadrature. In
particular, given the possibly singular nature of $W$ at corners,
\S\ref{sec:quadrature} motivates the use of geometrically graded
$h$-$p$ quadrature, because it is of high order for such integrands.
The approximation will be referred to as a discretisation.
Nomenclature used in the following discussion is shown in Figure
\ref{fig:nomen}.  (The distinction between \emph{mesh} and \emph{node}
(quadrature) points is made in \S\ref{sec:gmfq}.)

Consider the contour integral of an arbitrary integrand
$ g \left( \zeta \right) $ around $ \Gamma $:
\begin{eqnarray*}
  \oint_{\Gamma}
    g \left( \zeta \right)
  d\zeta.
\end{eqnarray*}
Parameterise $ \Gamma $ using $ \gamma : [0, 1] \mapsto \mathbb{C} $,
such that
$
  {\zeta}_0 \equiv \gamma \left( 0 \right)
  =
  \gamma \left( 1 \right)
  \equiv
  {\zeta}_N
$,
with argument $t$ increasing in an anticlockwise direction around
$\Gamma$. The contour integral is
\cite[page 168]{DavisRabinowitz:84}:%
\footnote{%
  This requires that $ \gamma $ is continuous, and that
  $\gamma|_{\left[ t_{i-1}, t_i \right]} $ is continuously
  differentiable for a finite partition
  $ 0 = t_0 < t_1 < \dots < t_n = 1 $.
}

\begin{eqnarray*}
  \oint_{\Gamma}
    g \left( \zeta \right)
  d\zeta
  =
  \int_{0}^{1}
    g \left( \gamma \left( t \right) \right)
    {\displaystyle \frac{\partial \gamma}{\partial t}} \left( t \right)
  dt.
\end{eqnarray*}

\vfill

\pagebreak

Approximate the integral using an $h$-$p$ quadrature rule
$ \left\lbrace t_j, w_j \right\rbrace $ with $N$ node points, defining
$ {\zeta}_j = \gamma \left( t_j \right) $ and
$
  {\dot{\gamma}}_j
  =
  {\displaystyle \frac{\partial \gamma}{\partial t}} \left( t_j \right)
$:
\begin{eqnarray*}
  \oint_{\Gamma}
    g \left( \zeta \right)
  d\zeta
  =
  \int_{0}^{1}
    g \left( \gamma \left( t \right) \right)
    {\displaystyle \frac{\partial \gamma}{\partial t}} \left( t \right)
  dt
  \approx
  \sum_{j=1}^{N}
    g \left( \zeta_j \right)
    {\dot{\gamma}}_j
    w_j.
\end{eqnarray*}
This formula can be applied to the Cauchy integrals. For some fixed
$z\in\Gamma$, consider the integrand
$
  g \left( \zeta \right)
  =
  \left( W \left( \zeta \right)
  -
  W \left( z \right) \right) / \left( \zeta - z \right)
$.
Let $ W_j = W \left( \zeta_j \right) $, and redefine
$w_j\Leftarrow {\dot{\gamma}}_j w_j$ to absorb $ {\dot{\gamma}}_j $.
The Cauchy integral is:
\begin{equation}
  \oint_{\Gamma}
    \displaystyle\frac
    {
      W \left( \zeta \right) - W \left( z \right)
    }{
      \zeta - z
    }
  d\zeta
  \approx
  \sum_{j=1}^{N}
    \displaystyle\frac
    {
      W_j - W \left( z \right)
    }{
      {\zeta}_j - z
    }
    w_j,
  \qquad
  z \in \Gamma.
  \label{eq:skunk}
\end{equation}
Approximation of the solution to the CBIE requires approximation of the
Cauchy integrals without using $ W \left( z \right) $. Instead,
$\hat{W} \left( z \right) \approx W \left( z \right) $ is constructed
from values of $W$ at the quadrature points. (\ref{eq:skunk}) is
discretised into a linear system of order $N$.

Begin by choosing a set of $N$ different values of $z$ from around the
boundary. These points are called the \emph{collocation} points.  It is
known from analysis in the case of uniform meshes that collocation
points must not lie on node points~\cite{ProssdorfRathsfeld:89}. The
natural choice is to take as the collocation points the $N$ midpoints
(in the sense of arc length) between the $N$ node points.%
\footnote{%
  The choice
  $ {\zeta}_{k - 1/2} = \left( {\zeta}_{k-1} + {\zeta}_k \right) / 2 $
  is explicitly \emph{not} used, as this assumes the contour is linear
  between points parameterised by $ t_{k-1} $ and $ t_k $.
}
Let $t_{k-1/2} = \left( t_{k-1} + t_k \right) / 2$,
${\zeta}_{k - 1/2} = \gamma \left( t_{k - 1/2} \right)$ and
$ W_{k - 1/2} = W \left( {\zeta}_{k - 1/2} \right)$,
for $ k = 1:N $. Interpolation from known values of $U$ at the
node and collocation points is used with the CBIE to approximate the
$ W_{k-1/2} $.

Define two complex $N$-vectors of $W$ at the node and collocation
points:
\begin{eqnarray*}
  \mathbf{W}
  =
  \left[
  \begin{array}{c}
    W_1 \\ \vdots \\ W_N
  \end{array}
  \right],
  \qquad
  {\mathbf{W}}'
  =
  \left[
  \begin{array}{c}
    W_{1/2} \\ \vdots \\ W_{N - 1/2}
  \end{array}
  \right].
\end{eqnarray*}
Also define the real $N$-vectors
$ \mathbf{U} = \Re \left( \mathbf{W} \right) $,
$ {\mathbf{U}}' = \Re \left( {\mathbf{W}}' \right) $ and
$ \mathbf{V} = \Im \left( \mathbf{W} \right) $.  From (\ref{eq:skunk}),
an order $N$ linear system for the components of $\mathbf{W}$ and
${\mathbf{W}}'$ is determined:
\begin{equation}
  \sum_{j=1}^N
    \displaystyle\frac
    {
      W_j - W_{k - 1/2}
    }{
      {\zeta}_j - {\zeta}_{k - 1/2}
    }
    w_j
  =
  0,
  \qquad
  k = 1:N.
  \label{eq:donkey}
\end{equation}

In summary, discretisation of Cauchy's integral formula leads to a
linear system, the solution to which is an approximation to $W$ at
the $N$ node points on $ \Gamma $. This approximation can be used to
approximate $W$, and hence $U$, at points within $ \Omega $.

\vfill

\pagebreak

%%%%%%%%%%%%%%%%%%%%%%%%%%%%%%%%%%%%%%%%%%%%%%%%%%%%%%%%%%%%%%%%%%%%%%%%

\subsubsection{Linear Interpolation of $W$ at the Collocation Points}
\label{sec:linintpoly}

If the $ W_j $ were known, (\ref{eq:donkey}) could be directly used to
interpolate the $ W_{k-1/2} $. However, although $ U_j $ and
$ U_{k-1/2} $ are known explicitly, $ V_j $ and $ V_{k-1/2} $ are
not. The CBIEM implicitly approximates the $ V_{k-1/2} $ by
interpolation from the as yet undetermined $ V_j $ at points
near $ {\zeta}_{k-1/2} $. That is, if a rule on $O$ points is
being used, choose $O$ terms from the sequence:
\begin{eqnarray*}
  \dots, \;
  {\zeta}_{k-3}, \;
  {\zeta}_{k-2}, \;
  {\zeta}_{k-1}; \;
  {\zeta}_k, \;
  {\zeta}_{k+1}, \;
  {\zeta}_{k+2}, \;
  \dots.
\end{eqnarray*}
The discretisation of the CBIE in (\ref{eq:donkey}), coupled with the $ 2N
$ knowns $ U_j $ and $ U_{k - 1/2} $, and the interpolation for the $
V_{k-1/2} $, allows the approximation of the $N$ unknowns $ V_j $.
(The $ V_{k - 1/2} $ are \emph{not} explicitly required to be
calculated.)

The simplest interpolation for the $ V_{k-1/2} $ is a linear one,
between $ {\zeta}_{k-1} $ and $ {\zeta}_k $, that is use:
\begin{equation}
  V_{k - 1/2}
  \approx
  \left(
    V_{k-1} + V_k
  \right)
  /
  2.
  \label{eq:cow}
\end{equation}
More sophisticated interpolations to the $ V_{k - 1/2} $ could involve
higher degree polynomials, splines, or trigonometric polynomials.
Initially, the linear choice will be used to illustrate the process. In
\S\ref{sec:higher}, the method is extended to higher degree polynomials.
For a particular problem ($ \Omega, f $), there is an optimal choice of
degree for the interpolation, as errors incurred by interpolation
increase with the degree, and eventually become of greater magnitude
than those due to discretisation.

The $ W_j $ are found by solving for their imaginary parts $ V_j $, and
combining these with the knowns $ U_j $.  The linear system is set up
as follows. Substituting $ W_j = U_j + i V_j $ and
$ W_{k - 1/2} = U_{k - 1/2} + i V_{k - 1/2} $ into (\ref{eq:donkey}),
using (\ref{eq:cow}), and collecting knowns and unknowns:
\begin{equation}
  \sum_{j=1}^N
    \displaystyle\frac
    {
      \frac{1}{2} V_{k-1}
      +
      \frac{1}{2} V_k
      -
      V_j
    }{
      {\zeta}_j - {\zeta}_{k - 1/2}
    }
    w_j
  =
  i
  \sum_{j=1}^N
    \displaystyle\frac
    {
      U_{k - 1/2}
      -
      U_j
    }{
      {\zeta}_j - {\zeta}_{k - 1/2}
    }
    w_j,
  \qquad
  k = 1:N.
  \label{eq:pig}
\end{equation}
This is a system of $N$ equations for the $N$ unknowns $ V_j $,
with RHS determined by the knowns $ U_j $ and $ U_{k-1/2} $ (and of
course the associated $ {\zeta}_j $ and $ {\zeta}_{k-1/2} $).

%%%%%%%%%%%%%%%%%%%%%%%%%%%%%%%%%%%%%%%%%%%%%%%%%%%%%%%%%%%%%%%%%%%%%%%%

\subsubsection{Solution of the Collocation Equations}
\label{sec:solcoll}

In order to write (\ref{eq:pig}) as a linear system, consider the LHS of
its $k$th equation:
\begin{equation}
  {\textstyle \frac{1}{2}}
  V_{k-1}
  \sum_{j=1}^N
    \displaystyle\frac
    {
      w_j
    }{
      {\zeta}_j - {\zeta}_{k - 1/2}
    }
  +
  {\textstyle \frac{1}{2}}
  V_k
  \sum_{j=1}^N
    \displaystyle\frac
    {
      w_j
    }{
      {\zeta}_j - {\zeta}_{k - 1/2}
    }
  -
  \sum_{j=1}^N
    \displaystyle\frac
    {
      V_j
      w_j
    }{
      {\zeta}_j - {\zeta}_{k - 1/2}
    }.
  \label{eq:ape}
\end{equation}
Define a matrix $ A \in \mathbb{C}^{N \times N}$:
\begin{eqnarray*}
  A
  =
  \left[
  \begin{array}{*{4}{c}}
    \displaystyle\frac{w_1}{{\zeta}_1 - \zeta_{1/2}} &
    \displaystyle\frac{w_2}{{\zeta}_2 - \zeta_{1/2}} &
    \dots &
    \displaystyle\frac{w_N}{{\zeta}_N - \zeta_{1/2}}
    \\
    \displaystyle\frac{w_1}{{\zeta}_1 - \zeta_{3/2}} &
    \displaystyle\frac{w_2}{{\zeta}_2 - \zeta_{3/2}} &
    \dots &
    \displaystyle\frac{w_N}{{\zeta}_N - \zeta_{3/2}}
    \\
    \vdots &
    \vdots &
    \ddots &
    \vdots
    \\
    \displaystyle\frac{w_1}{{\zeta}_1 - \zeta_{N - 1/2}} &
    \displaystyle\frac{w_2}{{\zeta}_2 - \zeta_{N - 1/2}} &
    \dots &
    \displaystyle\frac{w_N}{{\zeta}_N - \zeta_{N - 1/2}}
  \end{array}
  \right].
\end{eqnarray*}
Also define a set of $N$ scalars $ H_k $, for $ k = 1:N $ (the row
sums of $A$):
\begin{eqnarray*}
  H_k
  =
  \sum_{j=1}^N
    \displaystyle\frac
    {
      w_j
    }{
      {\zeta}_j - {\zeta}_{k - 1/2}
    }.
\end{eqnarray*}
The first two terms of (\ref{eq:ape}) are then
$ \frac{1}{2} V_{k-1} H_k + \frac{1}{2} V_k H_k $. As $ \Gamma $ is
closed, $ {\zeta}_0 \equiv {\zeta}_N $ and hence $ V_0 = V_N $. Define
$B \in \mathbb{C}^{N \times N}$:
\begin{equation}
  B
  =
  \displaystyle\frac{1}{2}
  \left[
  \begin{array}{*{6}{c}}
    H_1  &        &         &         & H_1 \\
    H_2  & H_2    &         &         &     \\
         & \ddots & \ddots  &         &      \\
         &        & H_{N-1} & H_{N-1} &      \\
         &        &         & H_N     & H_N
  \end{array}
  \right].
  \label{eq:chicken}
\end{equation}
Let $ C = B - A $, then the LHS of (\ref{eq:pig}) is $ C \mathbf{V} $.
Let $ {\bf 1}_N \in \mathbb{R}^{N}$ represent the real column vector
with all components unity, and the operation
$ \mathrm{diag} \left( \mathbf{x} \right) $ on $N$-vector
$\mathbf{x}$ create the diagonal matrix of order $N$ with the
diagonal entries being the respective components of $ \mathbf{x} $.
Defining
$
  \mathbf{d}
  =
  \mathrm{diag} \left( A {\bf 1}_N \right) {\mathbf{U}}'
  -
  A \mathbf{U}
$,
the system for $ \mathbf{V} $ is:
\begin{equation}
  C \mathbf{V}
  =
  i
  \mathbf{d}.
  \label{eq:owl}
\end{equation}
Attempting to directly solve the complex linear system in
(\ref{eq:owl}) fails, as $ \mathbf{V} $ is overdetermined in two
separate ways. Firstly, $ \mathbf{V} $ is purely real, that is
$\Im\left( \mathbf{V} \right) = \mathbf{0} $.  Partitioning
(\ref{eq:owl}) into real and imaginary components yields two purely
real linear systems, either one of which can be solved for what should
be the same solution $ \hat{\mathbf{V}} $. The system to be solved is:
$ \Re \left( C \right) \mathbf{V} = - \Im \left( \mathbf{d} \right) $
or $ \Im \left( C \right) \mathbf{V} = \Re \left( \mathbf{d} \right) $.
By redefining $ C \Leftarrow \Re \left( C \right) $ and
$\mathbf{d}\Leftarrow - \Im \left( \mathbf{d} \right) $, the first
choice for the solution of $ \mathbf{V} $ is made. The linear system is
thus:
\begin{equation}
  C
  \mathbf{V} = \mathbf{d}.
  \label{eq:sheep}
\end{equation}
The second way that (\ref{eq:owl}) is overdetermined is that
$\mathbf{V} $ is known only to within a constant.%
\footnote{%
  A scalar multiple of $ {\bf 1}_N $.
}
So, if direct solution of (\ref{eq:sheep}) is attempted, singularity
problems will occur, and the resultant $ \hat{\mathbf{V}} $ will be
infinite.%
\footnote{%
  Well, a numerical approximation to $ \infty $!
}
To accommodate this, arbitrarily%
\footnote{%
  For test problems, $ \hat{V}_N = V \left( {\zeta}_N \right) $
  is actually used.
}
set $ V_N = 0 $, and compute the rest of the components of $ \mathbf{V}
$ by subtracting rows in (\ref{eq:sheep}).  The result is a fully
determined order $ N - 1 $ linear system. The $N$ rows of
(\ref{eq:sheep}) are:
\begin{eqnarray*}
  {\left( C \mathbf{V} \right)}_j
  =
  {\mathbf{d}}_j,
  \qquad
  j = 1:N.
\end{eqnarray*}
Subtracting rows gives a system of $ N - 1 $ equations:
\begin{eqnarray*}
  {\left( C^{*} \mathbf{V} \right)}_{j-1}
  =
  {\left( C \mathbf{V} \right)}_j
  -
  {\left( C \mathbf{V} \right)}_{j-1}
  =
  {\mathbf{d}}_j - {\mathbf{d}}_{j-1}
  =
  {\mathbf{d}}_{j-1}^{*},
  \qquad
  j = 2:N.
\end{eqnarray*}
Defining $ {\mathbf{V}}^{*} \in \mathbb{R}^{N-1}$ as the first $N-1$
entries of $ \mathbf{V} $ (where the last entry is zero):
\begin{equation}
  C^{*}
  {\mathbf{V}}^{*}
  =
  {\mathbf{d}}^{*}.
  \label{eq:deer}
\end{equation}
Here:
\begin{eqnarray*}
  C^{*}
  & = &
  \left[
  \begin{array}{ccc}
    C_{2,1} - C_{1,1} & \dots & C_{2,N-1} - C_{1,N-1} \\
    C_{3,1} - C_{2,1} & \dots & C_{3,N-1} - C_{2,N-1} \\
    \vdots & \ddots & \vdots \\
    C_{N,1} - C_{N-1,1} & \dots & C_{N,N-1} - C_{N-1,N-1}
  \end{array}
  \right]
  \in
  \mathbb{R}^{(N-1) \times (N-1)}
  \\
  {\mathbf{d}}^{*}
  & = &
  \left[
  \begin{array}{c}
    d_2 - d_1 \\
    d_3 - d_2 \\
    \vdots \\
    d_N - d_{N-1}
  \end{array}
  \right]
  \in \mathbb{R}^{N-1}.
\end{eqnarray*}
Solution of the order $ N - 1 $ linear system in (\ref{eq:deer}) yields
the approximation to $ \mathbf{V} $, and hence $ \mathbf{W} $ ($W$ at
the $N$ node points).

%%%%%%%%%%%%%%%%%%%%%%%%%%%%%%%%%%%%%%%%%%%%%%%%%%%%%%%%%%%%%%%%%%%%%%%%

\subsection{Approximation of $U$ at Interior Points}
\label{sec:singsub}

The approximation $ \mathbf{W} = \mathbf{U} + i \mathbf{V} $ on
$\Gamma$ is used to approximate $W$ (and hence $U$) at interior points
of $\Omega$, using Cauchy's integral formula for points within
$\Omega$:
\begin{eqnarray*}
  W \left( z \right)
  =
  \displaystyle\frac{1}{2 \pi i}
  \oint_{\Gamma}
    \displaystyle\frac
    {
      W \left( \zeta \right)
    }{
      \zeta - z
    }
  d\zeta,
  \qquad
  z \in \Omega.
\end{eqnarray*}
Discretising this gives the approximation, for $ z \in \Omega $:
\begin{equation}
  W \left( z \right)
  \approx
  \displaystyle\frac{1}{2 \pi i}
  \sum_{j=1}^N
    \displaystyle\frac
    {
      W \left( {\zeta}_j \right)
    }{
      {\zeta}_j - z
    }
    w_j.
  \label{eq:egret}
\end{equation}
Numerical problems occur using this simple approximation, as points
$ z \in \Omega $ near the boundary (where $ z - {\zeta}_j $ is small,
for some $ {\zeta}_j $), generate very large terms in the sum.
In fact, the integrand is nearly singular, so only poor accuracy is
expected. Instead, the technique of singularity subtraction (referenced
in \S\ref{sec:cbie}) uses the result from (\ref{eq:skunk}):
\begin{eqnarray*}
  \displaystyle\frac{1}{2 \pi i}
  \sum_{j=1}^N
    \displaystyle\frac
    {
      W \left( {\zeta}_j \right)
      -
      W \left( z \right)
    }{
      {\zeta}_j - z
    }
    w_j
  =
  0.
\end{eqnarray*}
The integrand is now smooth, and good results can be expected from
quadrature. This yields $ W \left( z \right) $ as a `corrected'
(\ref{eq:egret}):
\begin{eqnarray*}
  W \left( z \right)
  =
  \left[
    \sum_{j=1}^N
      \displaystyle\frac
      {
        W \left( {\zeta}_j \right)
      }{
        {\zeta}_j - z
      }
      w_j
  \right]
  /
  \left[
    \sum_{j=1}^N
      \displaystyle\frac
      {
        1
      }{
        {\zeta}_j - z
      }
      w_j
  \right].
\end{eqnarray*}
Implementation of the CBIEM using this result is successful. The code
supplied (see \S\ref{sec:results}) does not go beyond the stage of the
computation of $ \hat{\mathbf{V}} $ on $ \Gamma $, as it is known that
the approximation to $U$ in the interior is actually \emph{more}
accurate than the approximations to $V$ on the
boundary~\cite{Chandler:90}.  Experiment demonstrates this, and thus
computation of $U$ at interior points need not be further described.

\pagebreak

%%%%%%%%%%%%%%%%%%%%%%%%%%%%%%%%%%%%%%%%%%%%%%%%%%%%%%%%%%%%%%%%%%%%%%%%

\subsection{Performance of the CBIEM on Model Problems}
\label{sec:typedir}

Application of the CBIEM yields different quality results depending on
the continuity of the model problem and the contour. In the simplest
case, both are smooth, there is no singularity, and standard quadrature
gives good results. In fact, because of periodicity, even the
trapezoidal rule on a uniform mesh gives very good results. There is no
need to use $h$-$p$ quadrature, but it will still work well.

Now consider the case where $ \Gamma $ is smooth, but $W$ has
singularities. For example, if $
W \left( \gamma \left( s \right) \right) \sim s^{1/2} $,
the integrand of the CBIE (\ref{eq:seal}) can have behaviour
$ s^{-1/2} $ near the corner.  However, using singularity subtraction
and $h$-$p$ quadrature, experiment demonstrates that both the
discretisation error and the final error in $ \hat{\mathbf{V}} $ on the
boundary are superb. (See also the example of complex contour
integration in \S\ref{sec:ccontint}.)

If $W$ is smooth, but $ \Gamma $ is not (has corners), in general,
the errors will be expected to increase with the sharpness of the
corner. The most difficult cases are cusps or reentrant corners (e.g.
the corner in a cardioid). Even for a model problem with a smooth
solution $U$, a corner singularity in $V$ (and hence $W$) will
occur.

Let $r$ represent radial distance from a corner. It is
known~\cite[pages 257--259]{StrangFix:73}, that at a corner with
interior angle $ \left( 1 - \chi \right) \pi $, a singularity of the
form $ r^{1/\left(1-\chi\right)} $ will be found. At worst, for a
reentrant corner, $ \chi = -1 $, so the form is $ r^{1/2} $.  A good
model problem is thus a contour with a corner where the true solution
has local behaviour $ U \sim r^{1/2} $.  This is obtained, for example,
using $ W \left( z \right) = z^{1/2} $.  If a uniform mesh were used,
the greatest component of the error in $V$ will come from the intervals
adjacent to the corner. The use of a geometrically graded mesh reduces
this component to a level comparable with that of other mesh
intervals.  Errors will not be of the very high order that is expected
for smooth contours, but should still be acceptable.

\pagebreak

%%%%%%%%%%%%%%%%%%%%%%%%%%%%%%%%%%%%%%%%%%%%%%%%%%%%%%%%%%%%%%%%%%%%%%%%

\subsection{Higher Degree Interpolatory Polynomials}

\subsubsection{Introduction}
\label{sec:higher}

To extend the technique described in \S\ref{sec:linintpoly}, the linear
interpolation in (\ref{eq:cow}) is replaced by a higher degree
interpolation. The net result of this is to change the definition of
the matrix $B$ in (\ref{eq:chicken}). Other possible techniques of
improving the accuracy of the interpolation, such as splines, are not
considered here, as they are difficult to implement. The principle
involved is that increasing the order of the interpolatory polynomial
should reduce the discretisation error, which is expected to be
greatest on the largest intervals.

The `nearest' $O$ node points on either side of $ {\zeta}_{k-1/2} $
are used to yield an interpolatory polynomial of degree $ O - 1 $.
In general, for points far from the nearest corner, this means to take
the first $O$ terms of the sequence
$ {\zeta}_k, {\zeta}_{k-1}, {\zeta}_{k+1}, {\zeta}_{k-2}, \dots $.
Otherwise, the term `nearest' is used loosely, as interpolation
cannot continue around a corner. Where there are less than $ O/2 $ node
points between the collocation point and the nearest corner (including
the node on the corner); instead $O$ points from and including
the corner are used. This results in an interpolatory polynomial that
is expected to be least accurate at the collocation point adjacent to
the corner. (A possible improvement in this schema is to
organise the interpolation rules such that their order increases say,
linearly, with node index away from the corner.) The constraint on
$O$ due to the mesh parameter $D$ is described in \S\ref{sec:limO}.

\S\S\ref{sec:lagform} and \ref{sec:Bconst} deal with technical
implementation issues, and can be skipped without loss of continuity.

%%%%%%%%%%%%%%%%%%%%%%%%%%%%%%%%%%%%%%%%%%%%%%%%%%%%%%%%%%%%%%%%%%%%%%%%

\subsubsection{Lagrange Form of the Interpolatory Polynomial}
\label{sec:lagform}

The Lagrange form of the interpolatory polynomial is used.
Given a set of values for the $ V_j $, at positions
$ {\zeta}_j = \gamma \left( t_j \right) $, the approximation to $V$ at
point $ {\zeta}_{k-1/2} $ is:
\begin{eqnarray*}
  V \left( {\zeta}_{k-1/2} \right)
  \approx
  \hat{V}_{k-1/2}
  =
  \sum_{j \in F_k}
    {\lambda}_j
    \left( t_{k-1/2} \right)
    V_j.
\end{eqnarray*}
Here:\footnote{Warning: Replacing $ t_{\nu} $ with $ {\zeta}_{\nu} $
and $ t_j $ with $ {\zeta}_j $ in this formula \emph{cannot} be done, as
the contour segments are \emph{not} necessarily straight.}
\begin{eqnarray*}
  {\lambda}_j
  \left(
    t
  \right)
  =
  \prod_{\nu \in F_k, \; \nu \ne j}
    \displaystyle\frac
    {
      t - t_{\nu}
    }{
      t_j - t_{\nu}
    }.
\end{eqnarray*}
$ F_k $ is a set of indices of the nearest node points, specific to
the collocation point $ {\zeta}_{k-1/2} $. Specifically, where
a degree $ O - 1 $ interpolatory polynomial is used on $O$
points, and $C$ is the index of the nearest corner to
$ {\zeta}_{k-1/2} $:
\begin{eqnarray*}
  F_k
  =
  \left\lbrace
  \begin{array}{ll}
    \left\lbrace
      k - O/2, \dots, k-1, k, \dots, k + O/2 - 1
    \right\rbrace
    &
    \mathrm{~in~general}
    \\
    \left\lbrace
      C, \dots, C + O - 1
    \right\rbrace
    &
    \mathrm{~if~} k - O/2 < C
    \\
    \left\lbrace
      C - O + 1, \dots, C
    \right\rbrace
    &
    \mathrm{~if~} k + O/2 > C.
  \end{array}
  \right.
\end{eqnarray*}
Let $ F_k $ be the $k$th row of a table $F$, and let the above
$ {\lambda}_j $ be the $j$th element of the $k$th row of another
table, $L$, of the associated weights.  The notations
$ F \left( k, j \right) $ and $ L \left( k, j \right) $ are used to
describe the set of $O$ nodal indices and weights associated with the
interpolation at point $ {\zeta}_{k-1/2} $, where $ k = 1:N $ and
$ j = 1:O $. The structure of $F$ becomes more complicated with
increasing $O$ and with increasing number of corners. Details of its
construction are not provided here, but can be read from the program
\texttt{cbiem.m} in Appendix \ref{app:cbiem}.

\vfill

\pagebreak

For some integer $D$, on each segment of $ \Gamma $ (with a corner at
each end), a mesh is constructed that has $ D-1 $ internal points
between each corner and a (wide) interval which spans the centre of the
segment. This results in a total of $ 2D $ mesh points, and $ 2D - 1 $
mesh intervals (see Figure \ref{fig:nomen}).

For example, if $ D = 3 $, then there are $ 4 $ interior and $ 2 $
corner mesh points on each segment, together with $ 4 $ extra interior
points, for a total of $ S = {\left( D + 1 \right)}^2 + 1 = 10 $.
Further, if there are $ NC = 3 $ corners, then the linear system has
order $ N = NC (S - 1) = 27 $. If $ O = 6 $ (quintic interpolation
about the `nearest' $ 6 $ points), then
$ F \in {\mathbb{N}}^{N \times O} $.  Where divisions in the structure
of $F$ due to the corners are reflected by partitions, $F$
is:%
\footnote{%
  Liberal use of \textsc{matlab} notation is made, and there
  is a confusion between computer array element and mathematical matrix
  entry notation: $ F \left( j, k \right) \equiv F_{j, k} $.
}

\tiny
\begin{eqnarray*}
  \mbox{\normalsize ~$ F = $~}
  \left[
  \begin{array}{*{6}{c}}
    27  &1  &2  &3  &4  &5  \\
    27  &1  &2  &3  &4  &5  \\
    27  &1  &2  &3  &4  &5  \\
    1  &2  &3  &4  &5  &6  \\
    2  &3  &4  &5  &6  &7  \\
    3  &4  &5  &6  &7  &8  \\
    4  &5  &6  &7  &8  &9  \\
    4  &5  &6  &7  &8  &9  \\
    4  &5  &6  &7  &8  &9  \\
    \hline
    9  &10  &11  &12  &13  &14  \\
    9  &10  &11  &12  &13  &14  \\
    9  &10  &11  &12  &13  &14  \\
    10  &11  &12  &13  &14  &15  \\
    11  &12  &13  &14  &15  &16  \\
    12  &13  &14  &15  &16  &17  \\
    13  &14  &15  &16  &17  &18  \\
    13  &14  &15  &16  &17  &18  \\
    13  &14  &15  &16  &17  &18  \\
    \hline
    18  &19  &20  &21  &22  &23  \\
    18  &19  &20  &21  &22  &23  \\
    18  &19  &20  &21  &22  &23  \\
    19  &20  &21  &22  &23  &24  \\
    20  &21  &22  &23  &24  &25  \\
    21  &22  &23  &24  &25  &26  \\
    22  &23  &24  &25  &26  &27  \\
    22  &23  &24  &25  &26  &27  \\
    22  &23  &24  &25  &26  &27
  \end{array}
  \right]
  \mbox{\normalsize ~$ = $~}
  \left[
  \begin{array}{c}
    \mbox{\normalsize $ Ft(1) $} \\
    \mbox{\normalsize $ Ft(2) + S *\mbox{ones}(Ft) $} \\
    \mbox{\normalsize $ \vdots $} \\
    \mbox{\normalsize $ Ft(NC) + (NC-1) S*\mbox{ones}(Ft) $}
  \end{array}
  \right].
\end{eqnarray*}
\normalsize

For $ j = 1:S-1 $ and $ k = 1:O $, in general
$ Ft \in \mathbb{N}^{\left( S-1 \right) \times O} $ is:
\begin{eqnarray*}
  Ft \left( j, k \right)
  =
  \left[
  \begin{array}{ll}
    k - 1    & j \leqslant O/2 \\
    k - 1 + j - O/2  & O/2 < j < S - O/2 \\
    k - 1 + S - O  & S - O/2 \leqslant j
  \end{array}
  \right].
\end{eqnarray*}
(Exception: $ F(1:O/2,1) \Leftarrow N \times {\bf 1}_{O/2} $.)

%%%%%%%%%%%%%%%%%%%%%%%%%%%%%%%%%%%%%%%%%%%%%%%%%%%%%%%%%%%%%%%%%%%%%%%%

\subsubsection{Construction of $B$}
\label{sec:Bconst}

The matrix $B$ is required in the construction of the linear system
in (\ref{eq:owl}), and is the only thing that changes when $O$ is
varied. For the case of linear interpolatory polynomials ($ O = 2 $),
a formula involving the terms $ H_k $ is used to approximate the value
of $V$ at the collocation points.  For larger $O$, this is replaced
with a considerably more sophisticated formula. As $O$
increases, the bandwidth of $B$ increases. Naturally, this new $B$
simplifies to the earlier definition if $ O = 2 $ is used, but is
obtained at greater computational expense.

The crucial change is in the approximation to the value of
$ V_{k-1/2} $, which in
(\ref{eq:pig}) is the term
$ \frac{1}{2} V_{k-1} + \frac{1}{2} V_k $. This is now replaced with:
\begin{eqnarray*}
  \sum_{j \in F_k}
    {\lambda}_j
    \left(
      {\zeta}_{k-1/2}
    \right)
    V_j
  =
  \sum_{j = 1}^O
    L \left( k, j \right)
    V_{F \left( k, j \right)}.
\end{eqnarray*}
The first $O$ terms of the LHS of the $k$th line of
(\ref{eq:pig}) (there is only one other term) are now:
\begin{eqnarray*}
  \sum_{j = 1}^O
    L \left( k, j \right)
    V_{F \left( k, j \right)}
  H_k.
\end{eqnarray*}
Construction of the real order $N$ matrix with $ \left( k, j \right)$th
entry  $ L \left( k, j \right) V_{F \left( k, j \right)} $ is required.
Multiplying each row by the corresponding $ H_k $ converts this to
$B$.

Details of the construction of $B$ are not provided here, but the
illustrative example used in \S\ref{sec:lagform} is continued. Recall
that $ D = 3 $, $ NC = 3 $ and $ O = 6 $, so $ S = 10 $, and $N=27$.
First construct $ Bt $ (a `skewed' version of $L$), and then calculate
$B$ by multiplying $ Bt $ through by the $ H_k $.  $ Bt $ is
constructed using a `shift vector' $G$, where $ G_k $ is equal to the
number of zeros to be put in front of row $k$ of $L$ to make it row $k$
of $ Bt $. This $G$ has a structure formed from a temporary $ Gt $:
\begin{eqnarray*}
  Gt
  & = &
  {\left[
    \mathrm{zeros} \left( 1, O/2 \right)
    \quad
    \left[ 1:S-O-1 \right]
    \quad
    \left( S - O \right) {\bf 1}_{O/2}^{\top}
  \right]}^{\top}
  \\
  G
  & = &
  {\left[
    Gt
    \quad
    Gt + S-1
    \quad
    \dots
    \quad
    Gt + NC \left( S-1 \right)
  \right]}^{\top}
  \\
  Bt \left( k, : \right)
  & = &
  {\left[
    \mathrm{zeros} \left( 1, G \left( k \right) \right)
    \; \;
    L \left( k, : \right)
    \; \;
    \mathrm{zeros} \left( 1, N - G \left( k \right) - O/2 \right)
  \right]}^{\top}
  \quad
  k = 1:N
  \\
  B
  &
  \Leftarrow
  &
  Bt \left( :, 2:N \right)
  \qquad
  B \left( 1:O/2, N \right)
  \;
  \Leftarrow
  \;
  Bt \left( 1:O/2, 1\right).
\end{eqnarray*}

\begin{figure}[htbp]
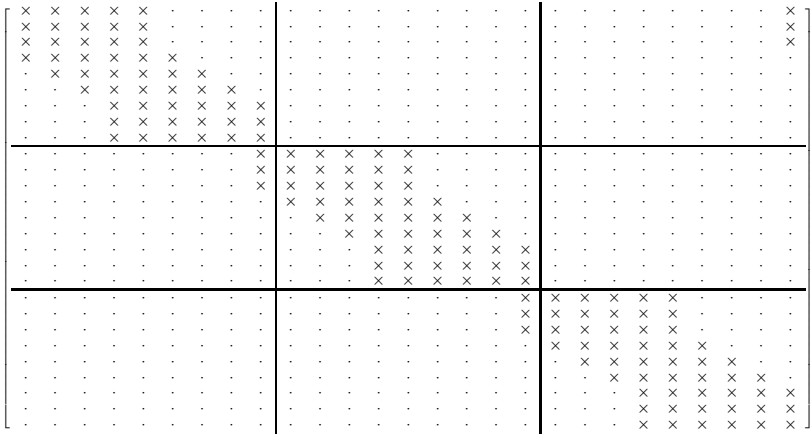

  \tiny
  \begin{eqnarray*}
  \left[
  \begin{array}{*{2}{*{9}{@{\hspace{1mm}}c@{\hspace{1mm}}}|}*{9}{@{\hspace{1mm}}c@{\hspace{1mm}}}}
    \times&\times&\times&\times&\times&\cdot&\cdot&\cdot&\cdot&\cdot&\cdot&\cdot&\cdot&\cdot&\cdot&\cdot&\cdot&\cdot&\cdot&\cdot&\cdot&\cdot&\cdot&\cdot&\cdot&\cdot&\times\\
    \times&\times&\times&\times&\times&\cdot&\cdot&\cdot&\cdot&\cdot&\cdot&\cdot&\cdot&\cdot&\cdot&\cdot&\cdot&\cdot&\cdot&\cdot&\cdot&\cdot&\cdot&\cdot&\cdot&\cdot&\times\\
    \times&\times&\times&\times&\times&\cdot&\cdot&\cdot&\cdot&\cdot&\cdot&\cdot&\cdot&\cdot&\cdot&\cdot&\cdot&\cdot&\cdot&\cdot&\cdot&\cdot&\cdot&\cdot&\cdot&\cdot&\times\\
    \times&\times&\times&\times&\times&\times&\cdot&\cdot&\cdot&\cdot&\cdot&\cdot&\cdot&\cdot&\cdot&\cdot&\cdot&\cdot&\cdot&\cdot&\cdot&\cdot&\cdot&\cdot&\cdot&\cdot&\cdot\\
    \cdot&\times&\times&\times&\times&\times&\times&\cdot&\cdot&\cdot&\cdot&\cdot&\cdot&\cdot&\cdot&\cdot&\cdot&\cdot&\cdot&\cdot&\cdot&\cdot&\cdot&\cdot&\cdot&\cdot&\cdot\\
    \cdot&\cdot&\times&\times&\times&\times&\times&\times&\cdot&\cdot&\cdot&\cdot&\cdot&\cdot&\cdot&\cdot&\cdot&\cdot&\cdot&\cdot&\cdot&\cdot&\cdot&\cdot&\cdot&\cdot&\cdot\\
    \cdot&\cdot&\cdot&\times&\times&\times&\times&\times&\times&\cdot&\cdot&\cdot&\cdot&\cdot&\cdot&\cdot&\cdot&\cdot&\cdot&\cdot&\cdot&\cdot&\cdot&\cdot&\cdot&\cdot&\cdot\\
    \cdot&\cdot&\cdot&\times&\times&\times&\times&\times&\times&\cdot&\cdot&\cdot&\cdot&\cdot&\cdot&\cdot&\cdot&\cdot&\cdot&\cdot&\cdot&\cdot&\cdot&\cdot&\cdot&\cdot&\cdot\\
    \cdot&\cdot&\cdot&\times&\times&\times&\times&\times&\times&\cdot&\cdot&\cdot&\cdot&\cdot&\cdot&\cdot&\cdot&\cdot&\cdot&\cdot&\cdot&\cdot&\cdot&\cdot&\cdot&\cdot&\cdot\\
    \hline
    \cdot&\cdot&\cdot&\cdot&\cdot&\cdot&\cdot&\cdot&\times&\times&\times&\times&\times&\times&\cdot&\cdot&\cdot&\cdot&\cdot&\cdot&\cdot&\cdot&\cdot&\cdot&\cdot&\cdot&\cdot\\
    \cdot&\cdot&\cdot&\cdot&\cdot&\cdot&\cdot&\cdot&\times&\times&\times&\times&\times&\times&\cdot&\cdot&\cdot&\cdot&\cdot&\cdot&\cdot&\cdot&\cdot&\cdot&\cdot&\cdot&\cdot\\
    \cdot&\cdot&\cdot&\cdot&\cdot&\cdot&\cdot&\cdot&\times&\times&\times&\times&\times&\times&\cdot&\cdot&\cdot&\cdot&\cdot&\cdot&\cdot&\cdot&\cdot&\cdot&\cdot&\cdot&\cdot\\
    \cdot&\cdot&\cdot&\cdot&\cdot&\cdot&\cdot&\cdot&\cdot&\times&\times&\times&\times&\times&\times&\cdot&\cdot&\cdot&\cdot&\cdot&\cdot&\cdot&\cdot&\cdot&\cdot&\cdot&\cdot\\
    \cdot&\cdot&\cdot&\cdot&\cdot&\cdot&\cdot&\cdot&\cdot&\cdot&\times&\times&\times&\times&\times&\times&\cdot&\cdot&\cdot&\cdot&\cdot&\cdot&\cdot&\cdot&\cdot&\cdot&\cdot\\
    \cdot&\cdot&\cdot&\cdot&\cdot&\cdot&\cdot&\cdot&\cdot&\cdot&\cdot&\times&\times&\times&\times&\times&\times&\cdot&\cdot&\cdot&\cdot&\cdot&\cdot&\cdot&\cdot&\cdot&\cdot\\
    \cdot&\cdot&\cdot&\cdot&\cdot&\cdot&\cdot&\cdot&\cdot&\cdot&\cdot&\cdot&\times&\times&\times&\times&\times&\times&\cdot&\cdot&\cdot&\cdot&\cdot&\cdot&\cdot&\cdot&\cdot\\
    \cdot&\cdot&\cdot&\cdot&\cdot&\cdot&\cdot&\cdot&\cdot&\cdot&\cdot&\cdot&\times&\times&\times&\times&\times&\times&\cdot&\cdot&\cdot&\cdot&\cdot&\cdot&\cdot&\cdot&\cdot\\
    \cdot&\cdot&\cdot&\cdot&\cdot&\cdot&\cdot&\cdot&\cdot&\cdot&\cdot&\cdot&\times&\times&\times&\times&\times&\times&\cdot&\cdot&\cdot&\cdot&\cdot&\cdot&\cdot&\cdot&\cdot\\
    \hline
    \cdot&\cdot&\cdot&\cdot&\cdot&\cdot&\cdot&\cdot&\cdot&\cdot&\cdot&\cdot&\cdot&\cdot&\cdot&\cdot&\cdot&\times&\times&\times&\times&\times&\times&\cdot&\cdot&\cdot&\cdot\\
    \cdot&\cdot&\cdot&\cdot&\cdot&\cdot&\cdot&\cdot&\cdot&\cdot&\cdot&\cdot&\cdot&\cdot&\cdot&\cdot&\cdot&\times&\times&\times&\times&\times&\times&\cdot&\cdot&\cdot&\cdot\\
    \cdot&\cdot&\cdot&\cdot&\cdot&\cdot&\cdot&\cdot&\cdot&\cdot&\cdot&\cdot&\cdot&\cdot&\cdot&\cdot&\cdot&\times&\times&\times&\times&\times&\times&\cdot&\cdot&\cdot&\cdot\\
    \cdot&\cdot&\cdot&\cdot&\cdot&\cdot&\cdot&\cdot&\cdot&\cdot&\cdot&\cdot&\cdot&\cdot&\cdot&\cdot&\cdot&\cdot&\times&\times&\times&\times&\times&\times&\cdot&\cdot&\cdot\\
    \cdot&\cdot&\cdot&\cdot&\cdot&\cdot&\cdot&\cdot&\cdot&\cdot&\cdot&\cdot&\cdot&\cdot&\cdot&\cdot&\cdot&\cdot&\cdot&\times&\times&\times&\times&\times&\times&\cdot&\cdot\\
    \cdot&\cdot&\cdot&\cdot&\cdot&\cdot&\cdot&\cdot&\cdot&\cdot&\cdot&\cdot&\cdot&\cdot&\cdot&\cdot&\cdot&\cdot&\cdot&\cdot&\times&\times&\times&\times&\times&\times&\cdot\\
    \cdot&\cdot&\cdot&\cdot&\cdot&\cdot&\cdot&\cdot&\cdot&\cdot&\cdot&\cdot&\cdot&\cdot&\cdot&\cdot&\cdot&\cdot&\cdot&\cdot&\cdot&\times&\times&\times&\times&\times&\times\\
    \cdot&\cdot&\cdot&\cdot&\cdot&\cdot&\cdot&\cdot&\cdot&\cdot&\cdot&\cdot&\cdot&\cdot&\cdot&\cdot&\cdot&\cdot&\cdot&\cdot&\cdot&\times&\times&\times&\times&\times&\times\\
    \cdot&\cdot&\cdot&\cdot&\cdot&\cdot&\cdot&\cdot&\cdot&\cdot&\cdot&\cdot&\cdot&\cdot&\cdot&\cdot&\cdot&\cdot&\cdot&\cdot&\cdot&\times&\times&\times&\times&\times&\times
    \end{array}
  \right]
  \end{eqnarray*}
  \normalsize
  \caption{Structure of $B$.}
  \label{fig:B}
\end{figure}

The overall structure of $B$ is depicted in Figure \ref{fig:B}, where
$\times$ and $\cdot$ represent nonzero and zero entries respectively.

The $ NC \times NC $ submatrices within the structure are each of order
$ S-1 $. The $k$th row generally consists of $ O = 6 $ contiguous
nonzero entries, starting at column $ G_k $.  These entries are the
$k$th row of $L$, multiplied by $ H_k $.  That is, the string of
$ O = 6 $ nonzero elements in row $k$ represents:

\begin{eqnarray*}
  \left[
    H_k L_{k,1} \quad
    H_k L_{k,2} \quad
    \dots \quad
    H_k L_{k,O-1} \quad
    H_k L_{k,O}
  \right].
\end{eqnarray*}

\pagebreak

Additionally, in the first $ O/2 $ rows, the first of these $O$
entries is shifted to the $N$th column, due to a `wraparound'
effect. In view of the banded structure of $B$, it appears that a
method designed to exploit this structure would be appropriate in the
solution of (\ref{eq:deer}). Unfortunately $A$ is neither banded, nor of
a particularly simple structure, so this approach is nontrivial, and
is a direction for further work.

%%%%%%%%%%%%%%%%%%%%%%%%%%%%%%%%%%%%%%%%%%%%%%%%%%%%%%%%%%%%%%%%%%%%%%%%

\subsubsection{Limitations on the Degree of the
     Interpolatory Polynomial}
\label{sec:limO}

The discussion and results presented in \S\ref{sec:quadrature} prompt
the use of a geometrically graded mesh, with the number of points in
the (closed) quadrature rule on each mesh interval linearly increasing
with interval number from the corner, beginning with $ 2 $ adjacent to
the corner, and becoming $D$ on the central (widest) interval. The
use of a linear grading is found in the literature of the finite
element method and the usual boundary element
method~\cite{BabuskaDorr:81,PostellStephan:90}.  Changing from linear
to quadratic or higher degree may reduce the errors, but its
implementation is beyond the scope of this report.

Consider a segment of $ \Gamma $ divided into a mesh on
$ 2 \left( D + 1 \right) $ points, including corners, with
$ m = 2 D + 1 $ intervals.  Basic quadrature rules on
$ n_j = 2, 3, \dots, D-1, D, D-1, \dots, 3, 2 $ points, are used over
intervals $ j = 1:m $.  After all of the common end points are
considered, the total number of points in the final composite
quadrature rule for that segment is $ S = {\left( D + 1 \right)}^2+1$.
Slightly different results would apply if the quadrature rules were
open. Recall that this choice is rejected, as it leads to node points
that avoid the singularity.

$S$ imposes a limit on $O$, the number of points used in the
interpolation rule. As $O$ is usually small (for reasons of
computational efficiency, typically $ O \leqslant 6 $), this limitation
is usually not significant. For example, if $ D = 2 $,
$ O \leqslant 8 $, and if $ D = 3 $, $ O \leqslant 16 $.

\pagebreak

%%%%%%%%%%%%%%%%%%%%%%%%%%%%%%%%%%%%%%%%%%%%%%%%%%%%%%%%%%%%%%%%%%%%%%%%

\section{Implementation and Results}
\label{sec:results}

\subsection{Implementation}

The CBIEM is implemented as a set of functions in \textsc{matlab}%
\footnote{%
  {\protect\textsc{matlab}} is an (interpreted) matrix
  computation package, and is a trademark of The Mathworks, Inc.
}
code,
presented in Appendix \ref{app:code}. (These functions appear in
alphabetic order, interspersed with several functions referenced in \S
\ref{sec:quadrature}.) The main routine is \texttt{cbiem.m}. The
parameterisation of the contour and its derivative are computed within
\texttt{cbiem.m}, and it calls an auxiliary function (\texttt{funccb.m})
to compute the true solution for a test problem.  The quadrature points
for basic rules are obtained by calling the function \texttt{gettw.m},
which provides either (closed) Newton--Cotes points (using a routine
internal to \texttt{gettw.m}), or calls another function,
\texttt{lobatto.m}, which computes the points for a closed
Gau{\ss}--Lobatto rule.  Two further functions are used to create an
$h$-$p$ composite quadrature rule out of a set of basic quadrature
rules (\texttt{hprmesh.m}), and compose a quadrature rule over a
contour with several corners (\texttt{rmesh.m}).  For testing
\texttt{cbiem.m} over a large set of parameters (e.g. generating the
data for Tables \ref{tab:res220} to \ref{tab:res2510}), a driver
function, \texttt{testcb.m}, is used.

Within \texttt{cbiem.m} is a description of its input parameters. For
test problems, where the true solution is known, it plots and
calculates norms of $ \mathbf{V} - \hat{\mathbf{V}} $, and also
calculates some discretisation errors.  Explicit computation of the
approximate solution within $ \Omega $ is \emph{not} performed.
Experiments with doing this demonstrate that the error results obtained
are of the same order as those returned.

\pagebreak

%%%%%%%%%%%%%%%%%%%%%%%%%%%%%%%%%%%%%%%%%%%%%%%%%%%%%%%%%%%%%%%%%%%%%%%%

\subsection{Experimental Results using a Teardrop Contour}

\subsubsection{Description}

\begin{figure}[htbp]
  \begin{center}
    \epsfig{file=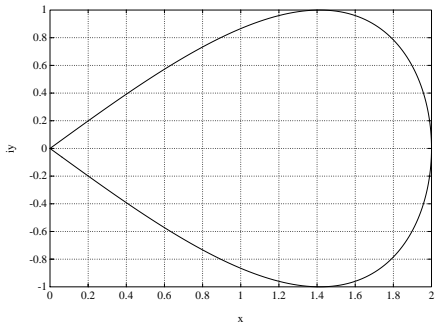,height=130pt}
    \caption{Teardrop contour used in the CBIEM experiments.}
    \label{fig:teardrop}
  \end{center}
\end{figure}

Although the CBIEM code is generalised to the situation of multiple
corners, good experimental contours have only one corner, to facilitate
isolation of the sources of error. This section describes numerical
results for the CBIEM, using a teardrop contour,\footnote{Another
important test contour is a cardioid with a reentrant corner.} depicted
in Figure \ref{fig:teardrop}.  The contour is parametrically given by:
\begin{eqnarray*}
  \gamma \left( t \right)
  =
  2 \sin \left( \pi t \right) + i \sin \left( 2 \pi t \right)
  \qquad
  t \in \left[ 0, 1 \right].
\end{eqnarray*}
This is the same contour as that used in \cite{Chandler:90}. It has a
right angle corner at the origin, which facilitates the use of test
problems $ W \left( z \right) = z^{\alpha} $. For
$ \alpha \in \left( 0, 1 \right) $, there is a discontinuity in the
derivative of the true solution at the origin, which becomes more
pathological as $ \alpha \to 0 $.

Error results are presented using an unweighted vector $ 2 $ norm:
\begin{eqnarray*}
  {|| \mathbf{V} - \hat{\mathbf{V}} ||}_2
  =
  {\textstyle \left[
    \sum_{i=1}^N
      {\left( V_i - \hat{V}_i \right)}^2
  \right]}^{1/2}.
\end{eqnarray*}
An appropriately weighted discretisation of the $ L_2 $ norm might seem
more appropriate, but would effectively only present the norm over the
central interval, as the widths of the end point intervals are very
small. The use of an infinity norm is also appealing, but the $ 2 $
norm allows the user to experiment with interpolation formula
gradings, to independently reduce the error over different regions of
the contour (see \S\ref{sec:improve}). Also, experimental data shows the
behaviour of the infinity norm is very similar to that of the $ 2 $
norm.

Tables \ref{tab:res220} to \ref{tab:res2510} present error results for
the teardrop contour, for three different model problems:
$ W \left( z \right) = z^2 $, $ z^{1/2} $ and $ z^{1/4} $; various
choices of two mesh grading parameters ($ \sigma $ and $D$); and
choices of $O$, the number of points used by the interpolatory
polynomial in the collocation process.  In each table
$ N = {\left( D + 1 \right)}^2 $ is the size of the linear system being
solved, such that there are $ 2 D + 1 $ mesh intervals between one
corner and the next (see \S\ref{sec:lagform}).  The nine tables cover
three illustrative choices of the mesh parameter $ \sigma $ for each of
the three model problems.  In each case, the results presented are for
a choice of $ \sigma $ close to the optimal $ \sigma $, and two nearby
values of $ \sigma $ that demonstrate the increase in the error in each
direction.  Results have been selected from a much larger data set.
Within each table, the minimum error result is emboldened.

\pagebreak

%%%%%%%%%%%%%%%%%%%%%%%%%%%%%%%%%%%%%%%%%%%%%%%%%%%%%%%%%%%%%%%%%%%%%%%%

\subsubsection{Observations of $ z^2 $}

This test function does \emph{not} have a singularity, and the results
are good.  Despite the corner, the error reduces with increasing either
$D$ or $O$, until a point is reached where roundoff error, caused
by excessive order in the interpolatory polynomial, begins to
encroach.

%%%%%%%%%%%%%%%%%%%%%%%%%%%%%%%%%%%%%%%%%%%%%%%%%%%%%%%%%%%%%%%%%%%%%%%%

% The first tables are for the nonsingular case z^2.

\begin{table}[htbp]
  \centering
  \begin{tabular}{||c||*{8}{c|}|}
    \hline\hline
    &\multicolumn{8}{c||}{$O$} \\
    \hline
    $N$ &8&10&12&14&16&18&20&22 \\
    \hline\hline
    25 & 4.6e-04 & 5.6e-05 & 6.0e-06 & 1.0e-06 &
      1.2e-06 & 1.2e-06 & 3.7e-06 & 3.0e-06 \\
    36 & 1.4e-04 & 1.3e-05 & 1.1e-06 & 8.0e-08 &
      4.2e-08 & 2.8e-08 & 7.8e-07 & 4.1e-07 \\
    49 & 7.1e-05 & 5.7e-06 & 3.9e-07 & 2.5e-08 &
      5.6e-09 & {\bf 3.5e-10} & 5.8e-08 & 2.9e-05 \\
    64 & 3.1e-05 & 2.1e-06 & 1.2e-07 & 6.5e-09 &
      1.5e-09 & 9.0e-09 & 2.6e-07 & 5.2e-02 \\
    81 & 1.6e-05 & 1.0e-06 & 5.5e-08 & 2.5e-09 &
      3.6e-10 & 1.2e-08 & 3.6e-05 & 2.0e-02 \\
    100& 8.6e-06 & 4.9e-07 & 2.3e-08 & 9.6e-10 &
      8.4e-10 & 1.0e-08 & 4.2e-05 & 1.6e-01 \\
    \hline\hline
  \end{tabular}
  \caption{
      $ \sigma = 0.20 $,
      $ U \left( z \right) = \Re \left( z^2 \right) $.
    }
  \label{tab:res220}
\end{table}

%%%%%%%%%%%%%%%%%%%%%%%%%%%%%%%%%%%%%%%%%%%%%%%%%%%%%%%%%%%%%%%%%%%%%%%%

\begin{table}[htbp]
  \centering
  \begin{tabular}{||c||*{8}{c|}|}
    \hline\hline
    &\multicolumn{8}{c||}{$O$} \\
    \hline
    $N$ &8&10&12&14&16&18&20&22 \\
    \hline\hline
    25 & 5.9e-05 & 6.6e-06 & 9.0e-06 & 9.3e-06 &
      9.3e-06 & 9.4e-06 & 9.3e-06 & 9.7e-06 \\
    36 & 1.5e-05 & 9.3e-07 & 1.8e-07 & 2.2e-07 &
      2.3e-07 & 2.3e-07 & 2.4e-07 & 1.4e-06 \\
    49 & 6.2e-06 & 3.1e-07 & 1.2e-08 & 5.0e-09 &
      5.5e-09 & 5.5e-09 & 6.9e-09 & 6.7e-08 \\
    64 & 2.6e-06 & 1.0e-07 & 4.0e-09 & 9.9e-11 &
      1.2e-10 & 1.3e-10 & 1.5e-10 & 4.7e-09 \\
    81 & 1.3e-06 & 4.5e-08 & 1.4e-09 & 4.1e-11 &
      2.3e-12 & 3.1e-12 & 4.6e-12 & 8.1e-09 \\
    100& 6.7e-07 & 2.0e-08 & 5.7e-10 & 1.4e-11 &
      {\bf 3.8e-13} & 5.2e-13 & 2.5e-11 & 4.6e-08 \\
    \hline\hline
  \end{tabular}
  \caption{
      $ \sigma = 0.28 $,
      $ U \left( z \right) = \Re \left( z^2 \right) $.
    }
  \label{tab:res228}
\end{table}

%%%%%%%%%%%%%%%%%%%%%%%%%%%%%%%%%%%%%%%%%%%%%%%%%%%%%%%%%%%%%%%%%%%%%%%%

\begin{table}[htbp]
  \centering
  \begin{tabular}{||c||*{8}{c|}|}
    \hline\hline
    &\multicolumn{8}{c||}{$O$} \\
    \hline
    $N$ &8&10&12&14&16&18&20&22 \\
    \hline\hline
    25 & 1.3e-04 & 1.4e-04 & 1.4e-04 & 1.4e-04 &
      1.4e-04 & 1.4e-04 & 1.4e-04 & 1.4e-04 \\
    36 & 6.3e-06 & 6.7e-06 & 6.8e-06 & 6.8e-06 &
      6.8e-06 & 6.8e-06 & 6.8e-06 & 6.8e-06 \\
    49 & 3.4e-07 & 3.1e-07 & 3.2e-07 & 3.2e-07 &
      3.2e-07 & 3.2e-07 & 3.2e-07 & 3.2e-07 \\
    64 & 1.1e-07 & 1.3e-08 & 1.5e-08 & 1.5e-08 &
      1.5e-08 & 1.5e-08 & 1.5e-08 & 1.5e-08 \\
    81 & 5.4e-08 & 7.9e-10 & 6.8e-10 & 6.9e-10 &
      6.9e-10 & 6.9e-10 & 6.9e-10 & 6.9e-10 \\
    100& 2.6e-08 & 3.7e-10 & {\bf 2.8e-11} & 3.1e-11 &
      3.1e-11 & 3.1e-11 & 3.1e-11 & 3.1e-11 \\
    \hline\hline
  \end{tabular}
  \caption{
      $ \sigma = 0.35 $,
      $ U \left( z \right) = \Re \left( z^2 \right) $.
    }
  \label{tab:res235}
\end{table}

\pagebreak

%%%%%%%%%%%%%%%%%%%%%%%%%%%%%%%%%%%%%%%%%%%%%%%%%%%%%%%%%%%%%%%%%%%%%%%%

\subsubsection{Observations of $ z^{1/2} $}

\begin{figure}[htbp]
  \begin{center}
    \epsfig{file=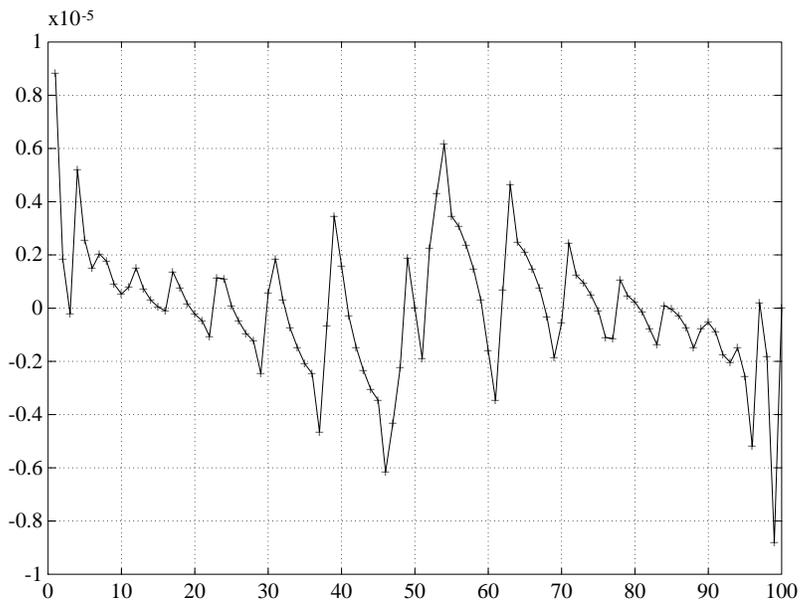,height=250pt}
    \caption{
      $ \protect\mathbf{V} - \hat{\protect\mathbf{V}} $ for
      \protect Table \ref{tab:res510}, using $ D = 9 $, $ O = 6 $.
    }
    \label{fig:vvsvhat}
  \end{center}
\end{figure}

This case has a corner singularity, and represents the `worst' that
singularities get in practice (that is, for $ z^{\alpha} $
singularities, in practice $ \alpha \geqslant 1/2 $). Although the
error decreases with increasing $D$ or $O$, it does so more slowly
than for $ z^2 $, and is orders of magnitude larger. As for $ z^2 $,
there comes a point where increasing $O$ causes the error to
increase, and indeed grow exponentially.  The results for
$\mathbf{V} - \hat{\mathbf{V}} $ for the case $ \sigma = 0.10 $,
$ D = 9 $, $ O = 6 $ are plotted in Figure \ref{fig:vvsvhat}. The
abscissae are plotted uniformly, for if they were plotted versus
parameter $t$, the geometric grading would bunch up most of the
results at the ends (corner of teardrop).  Observe that these results
are, as expected, antisymmetric.

% The second tables are for the singular case z^{1/2}.

\begin{table}[htbp]
  \centering
  \begin{tabular}{||c||*{5}{c|}|}
    \hline\hline
    &\multicolumn{5}{c||}{$O$} \\
    \hline
    $N$ &2&4&6&8&10 \\
    \hline\hline
  16 & 2.2e-02 & 1.1e-02 & 2.6e-02 & 7.3e-01 & 2.3e+00 \\
  25 & 1.4e-02 & 6.6e-03 & 8.5e-03 & 1.8e-01 & 4.5e+00 \\
  36 & 1.0e-02 & 2.2e-03 & 2.1e-03 & 4.3e-02 & 5.3e+00 \\
  49 & 8.1e-03 & 1.3e-03 & 9.3e-04 & 9.7e-03 & 1.8e+00 \\
  64 & 6.3e-03 & 6.8e-04 & 3.6e-04 & 2.2e-03 & 4.5e-01 \\
  81 & 4.9e-03 & 4.5e-04 & 1.9e-04 & 5.5e-04 & 1.0e-01 \\
  100& 3.9e-03 & 2.9e-04 & {\bf 1.1e-04} & 1.7e-04 & 2.2e-02 \\
    \hline\hline
  \end{tabular}
  \caption{
      $ \sigma = 0.05 $,
      $ U \left( z \right) = \Re \left( z^{1/2} \right) $.
    }
  \label{tab:res505}
\end{table}

%%%%%%%%%%%%%%%%%%%%%%%%%%%%%%%%%%%%%%%%%%%%%%%%%%%%%%%%%%%%%%%%%%%%%%%%

\begin{table}[htbp]
  \centering
  \begin{tabular}{||c||*{5}{c|}|}
    \hline\hline
    &\multicolumn{5}{c||}{$O$} \\
    \hline
    $N$ &2&4&6&8&10 \\
    \hline\hline
  16 & 2.0e-02 & 1.4e-02 & 1.3e-02 & 4.4e-02 & 2.7e-01 \\
  25 & 1.2e-02 & 5.2e-03 & 5.1e-03 & 1.4e-02 & 1.4e-01 \\
  36 & 8.3e-03 & 1.6e-03 & 1.5e-03 & 4.5e-03 & 4.9e-02 \\
  49 & 5.9e-03 & 6.2e-04 & 5.2e-04 & 1.4e-03 & 1.5e-02 \\
  64 & 4.5e-03 & 2.5e-04 & 1.6e-04 & 4.5e-04 & 4.9e-03 \\
  81 & 3.4e-03 & 1.3e-04 & 6.0e-05 & 1.4e-04 & 1.5e-03 \\
  100& 2.7e-03 & 8.4e-05 & {\bf 2.4e-05} & 4.6e-05 & 4.9e-04 \\
    \hline\hline
  \end{tabular}
  \caption{
      $ \sigma = 0.10 $,
      $ U \left( z \right) = \Re \left( z^{1/2} \right) $.
    }
  \label{tab:res510}
\end{table}

%%%%%%%%%%%%%%%%%%%%%%%%%%%%%%%%%%%%%%%%%%%%%%%%%%%%%%%%%%%%%%%%%%%%%%%%

\begin{table}[htbp]
  \centering
  \begin{tabular}{||c||*{5}{c|}|}
    \hline\hline
    &\multicolumn{5}{c||}{$O$} \\
    \hline
    $N$ &2&4&6&8&10 \\
    \hline\hline
  16 & 2.4e-02 & 2.3e-02 & 2.2e-02 & 2.8e-02 & 3.0e-02 \\
  25 & 1.2e-02 & 9.7e-03 & 9.3e-03 & 1.1e-02 & 1.3e-02 \\
  36 & 7.0e-03 & 3.7e-03 & 3.5e-03 & 4.3e-03 & 5.4e-03 \\
  49 & 4.7e-03 & 1.4e-03 & 1.3e-03 & 1.6e-03 & 2.1e-03 \\
  64 & 3.4e-03 & 5.7e-04 & 5.3e-04 & 6.5e-04 & 8.2e-04 \\
  81 & 2.6e-03 & 2.3e-04 & 2.0e-04 & 2.5e-04 & 3.1e-04 \\
  100& 2.0e-03 & 9.8e-05 & {\bf 8.0e-05} & 9.8e-05 & 1.2e-04 \\
    \hline\hline
  \end{tabular}
  \caption{
      $ \sigma = 0.15 $,
      $ U \left( z \right) = \Re \left( z^{1/2} \right) $.
    }
  \label{tab:res515}
\end{table}

\clearpage

%%%%%%%%%%%%%%%%%%%%%%%%%%%%%%%%%%%%%%%%%%%%%%%%%%%%%%%%%%%%%%%%%%%%%%%%

\subsubsection{Observations of $ z^{1/4} $}

A $ z^{1/4} $ singularity is beyond the range of singularities expected
for smooth test functions. The errors are worse again than for
$ z^{1/2} $, and they do not decrease as fast with increasing $D$ or
$O$. In fact, when the test problem is this pathological, the
Dirichlet problem is fast becoming a boundary layer problem, which
should be dealt with using more specialised methods.

% Thirdly, tables for z^{1/4}

\begin{table}[htbp]
  \centering
  \begin{tabular}{||c||*{4}{c|}|}
    \hline\hline
    &\multicolumn{4}{c||}{$O$} \\
    \hline
    $N$ &2&4&6&8 \\
    \hline\hline
    16 & 4.6e-02 & 4.4e-02 & 8.1e-01 & 1.0e+01 \\
    25 & 3.2e-02 & 3.3e-02 & 2.8e-01 & 3.2e+01 \\
    36 & 1.7e-02 & 1.9e-02 & 1.0e-01 & 4.9e+01 \\
    49 & 1.0e-02 & 1.2e-02 & 4.1e-02 & 3.3e+01 \\
    64 & 5.9e-03 & 6.9e-03 & 1.6e-02 & 1.3e+01 \\
    81 & 3.9e-03 & 4.2e-03 & 7.5e-03 & 1.8e+01 \\
    100& 2.9e-03 & {\bf 2.5e-03} & 3.6e-02 & 2.8e+02 \\
    \hline\hline
  \end{tabular}
  \caption{
    $ \sigma = 0.02 $,
    $ U \left( z \right) = \Re \left( z^{1/4} \right) $.
  }
  \label{tab:res2502}
\end{table}

%%%%%%%%%%%%%%%%%%%%%%%%%%%%%%%%%%%%%%%%%%%%%%%%%%%%%%%%%%%%%%%%%%%%%%%%

\begin{table}[htbp]
  \centering
  \begin{tabular}{||c||*{4}{c|}|}
    \hline\hline
    &\multicolumn{4}{c||}{$O$} \\
    \hline
    $N$ &2&4&6&8 \\
    \hline\hline
    16 & 3.6e-02 & 3.8e-02 & 5.8e-02 & 1.4e+00 \\
    25 & 2.0e-02 & 2.1e-02 & 3.0e-02 & 7.5e-01 \\
    36 & 1.0e-02 & 9.8e-03 & 1.4e-02 & 3.6e-01 \\
    49 & 5.4e-03 & 4.9e-03 & 7.1e-03 & 1.7e-01 \\
    64 & 3.2e-03 & 2.3e-03 & 3.3e-03 & 8.3e-02 \\
    81 & 2.1e-03 & 1.1e-03 & 1.6e-03 & 4.0e-02 \\
    100& 1.6e-03 & {\bf 5.6e-04} & 7.8e-04 & 2.0e-02 \\
    \hline\hline
  \end{tabular}
  \caption{
      $ \sigma = 0.05 $,
      $ U \left( z \right) = \Re \left( z^{1/4} \right) $.
    }
  \label{tab:res2505}
\end{table}

%%%%%%%%%%%%%%%%%%%%%%%%%%%%%%%%%%%%%%%%%%%%%%%%%%%%%%%%%%%%%%%%%%%%%%%%

\begin{table}[htbp]
  \centering
  \begin{tabular}{||c||*{4}{c|}|}
    \hline\hline
    &\multicolumn{4}{c||}{$O$} \\
    \hline
    $N$ &2&4&6&8 \\
    \hline\hline
    16 & 5.8e-02 & 6.0e-02 & 5.8e-02 & 6.4e-02 \\
    25 & 3.4e-02 & 3.5e-02 & 3.4e-02 & 3.7e-02 \\
    36 & 1.9e-02 & 1.9e-02 & 1.9e-02 & 2.0e-02 \\
    49 & 1.1e-02 & 1.1e-02 & 1.0e-02 & 1.1e-02 \\
    64 & 6.2e-03 & 6.2e-03 & 6.0e-03 & 6.5e-03 \\
    81 & 3.6e-03 & 3.5e-03 & 3.3e-03 & 3.7e-03 \\
    100& 2.1e-03 & 1.9e-03 & {\bf 1.8e-03} & 2.0e-03 \\
    \hline\hline
  \end{tabular}
  \caption{
      $ \sigma = 0.10 $,
      $ U \left( z \right) = \Re \left( z^{1/4} \right) $.
    }
  \label{tab:res2510}
\end{table}

\clearpage

%%%%%%%%%%%%%%%%%%%%%%%%%%%%%%%%%%%%%%%%%%%%%%%%%%%%%%%%%%%%%%%%%%%%%%%%

\subsubsection{The Black Art of Choosing $ \sigma $}

\begin{figure}[htbp]
  \begin{center}
    \epsfig{file=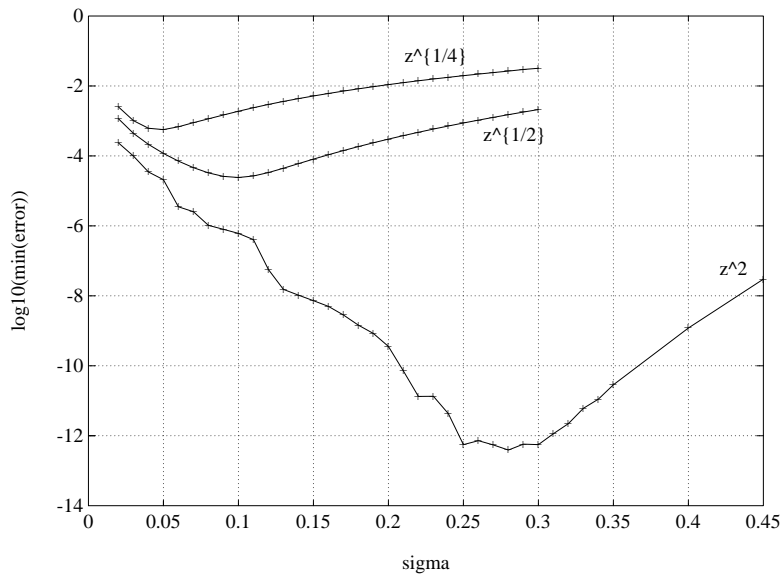,height=230pt}
    \caption{Variation of minimum error with $ \sigma $.}
    \label{fig:sigma}
  \end{center}
\end{figure}

The minimum error results obtained for each test problem are plotted
versus $ \sigma $ in Figure \ref{fig:sigma}, and demonstrate that there
is an optimal choice of $ \sigma $, which varies significantly with the
test problem.  The use of the CBIEM in applications, where the true
solution is not known in advance, could falter on the setting of $
\sigma $. If the computational cost is to be minimised, then it is
important to find the optimal $ \sigma $, however, it may be expensive
to try many $ \sigma $ until the optimal one is found.  The literature
does not justify a choice of $ \sigma $, but merely states it, e.g.
\cite{Stephan:88} uses $ \sigma = 0.15 $ for a particular (finite
element) application.  The optimal choice of $ \sigma $ for the
paradigm test problem $ z^{1/2} $ is $ \sigma = 0.10 $. As $ z^{1/2} $
is the worst singularity expected in practice (see
\S\ref{sec:typedir}), this should be a good guide as a starting guess
for any problem with an unknown solution.

%%%%%%%%%%%%%%%%%%%%%%%%%%%%%%%%%%%%%%%%%%%%%%%%%%%%%%%%%%%%%%%%%%%%%%%%

\subsubsection{Improvements in the Technique}
\label{sec:improve}

Consider Table \ref{tab:res510}, where the best error result is
obtained using $ O = 6 $. The error may be able to be reduced by
grading the order of the interpolatory polynomial over the mesh
intervals. Near the corner, the use of high order interpolation may
actually \emph{increase} the component of the error, although this may
be appropriate far away from the corner.  It may be ideal to grade the
order of the interpolatory polynomial from $ O = 2 $ near the corner,
to $ O = 6 $ (or greater) farthest from the corner.

\pagebreak

Direct implementation of this result requires extensive modification to
the matrix $B$ used by the CBIEM\footnote{These comments also refer
to the associated $F$ and $L$ matrices.} (see \S\ref{sec:higher}),
and is beyond the scope of this report.  Another way of achieving the
same effect is to calculate $B$ matrices $ B_2, B_4, \dots, B_{16} $
for $ O = 2, 4, \dots, 16 $, then construct a new $B$ from
appropriate rows of them, and insert this new $B$ at the relevant
point in the CBIEM.

\begin{figure}[htbp]
  \begin{center}
    \epsfig{file=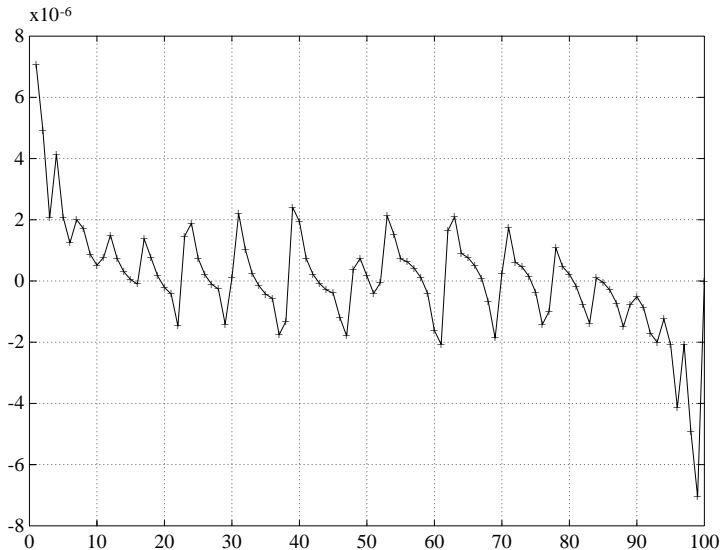,height=230pt}
    \caption{
      The best error results for $ z^{1/2} $, $ D = 9 $,
      $ \sigma = 0.10 $ and a graded interpolation rule.
      c.f. \protect Figure \ref{fig:vvsvhat}.
    }
    \label{fig:best}
  \end{center}
\end{figure}

However, intuition is misleading here. The minimum error in Table
\ref{tab:res510} is $ 2.4158 \times 10^{-5} $, using $ \sigma = 0.10 $
and $ D = 9 $. This corresponds to $ O = 6 $ on each mesh interval.
Many experiments in variation of the order of the interpolation rule,
holding fixed $ \sigma $ and $D$, find that the very best error
result that can be obtained is $ 1.7444 \times 10^{-5} $ (a $ 28\% $
reduction), using a grading with interpolation rules of
$ O = 12, 2, 2, 6, 6, 6, 10, 10, 10 $ over the $ 9 $ mesh intervals
from the corner to the centre. Surprisingly, the component of the error
over the first interval \emph{decreases} with increasing $O$.  This
is depicted in Figure \ref{fig:best}, which shows that the error is
uniformly distributed around the contour, except for the largest
component, at the corner.

It appears that what is happening is that the method has come up
against a discretisation error barrier. For this problem, the
discretisation error does not decrease particularly quickly, and is a
maximum at the corner.

\vfill

\pagebreak

%%%%%%%%%%%%%%%%%%%%%%%%%%%%%%%%%%%%%%%%%%%%%%%%%%%%%%%%%%%%%%%%%%%%%%%%

\section{Further Directions for Research}
\label{sec:further}

This section enumerates various possibilities for future work on the
CBIEM.

\begin{enumerate}
\item
  There is the potential for error reduction using graded
  interpolation rules (see \S\ref{sec:improve}). Similarly, other
  choices for the grading of the quadrature rules may assist in
  error reduction, e.g.  quadratic increase in degree of
  quadrature rule with node number from the corner, rather than
  linear as is presently used.
\item
  A proper examination of the computational efficiency of the
  CBIEM is required. This would involve setting up, say, a
  finite difference solution for the Dirichlet problem, and
  comparing flop counts required to obtain comparable accuracies.
\item
  Analysis of the choice of optimal $ \sigma $ is desirable.
  Currently, the method is hampered by this not being known in
  advance.
\item
  Application of the technique to conformal
  mapping~\cite{DeLilloElcratt:91} may be worthwhile.
\item
  It would be computationally efficient if solution of the linear
  system involving the matrix $ C = B - A $ could exploit the
  banded structure of $B$ (see \S\ref{sec:Bconst}).
\item
  An alternative collocation technique is
  possible~\cite{ProssdorfRathsfeld:89}.  Given nodes
  $ {\zeta}_j $ with weights $ w_j $, and collocation points
  $ {\zeta}_{j+1/2} $ with weights $ w_{j+1/2} $:
  \begin{eqnarray*}
    \oint_{\Gamma}
      F \left( \zeta \right)
    d\zeta
    \approx
    \sum_{j=1}^N
      F \left( {\zeta}_j \right)
      w_j
    \approx
    \sum_{k=1}^N
      F \left( {\zeta}_{k - 1/2} \right)
      w_{k-1/2}.
  \end{eqnarray*}
  Approximate the unknowns $ V_j = V \left( {\zeta}_j \right) $ by
  collocating at $ {\zeta}_{k-1/2} $, and the unknowns
  $ V_{k-1/2} = V \left( {\zeta}_{k-1/2} \right) $ by collocating at
  $ {\zeta}_j $. This gives an order $ 2N $ system for the $ 2N $
  unknowns $ V_j $ and $ V_{k-1/2} $, but avoids interpolation.
  \begin{eqnarray*}
    0
    & = &
    \sum_{j=1}^N
      \displaystyle\frac
      {
        W_j - W_{k-1/2}
      }{
        {\zeta}_j - {\zeta}_{k-1/2}
      }
      w_j
    \\
    0
    & = &
    \sum_{k=1}^N
      \displaystyle\frac
      {
        W_{k-1/2} - W_j
      }{
        {\zeta}_{k-1/2} - {\zeta}_j
      }
      w_{k-1/2}.
  \end{eqnarray*}
  The technique appears to be computationally wasteful, but may
  be worth investigating, as it would be simpler to implement.

\pagebreak

\item
  The CVBEM was developed to solve $ 2 $D fluid flow
  problems~\cite{Hromadka:84},%
  \footnote{%
    Other references to the CVBEM include
    \cite{%
      Gartland:88,%
      HromadkaGuymon:84a,%
      HromadkaGuymon:84b,%
      SchultzHong:89%
  }.}
  where components of the complex potential (the fluid potential
  $ \Phi $ or the streamline function $ \Psi $) are known at
  different points around the contour, typically from physical
  measurements.%
  \footnote{%
    Warning: the notation used here
    $ W = U + i V $ is equivalent to the fluid flow notation
    $ W = \Phi + i \Psi $, so that $U$ and $V$ here have a
    different meaning from the fluid flow case, where they are
    commonly the components of the velocity
    $ q = U \hat{\mathbf{i}} + V \hat{\mathbf{j}} $, and
    $
      U
      =
      {\displaystyle \frac{\partial \Phi}{\partial x}}
      =
      {\displaystyle \frac{\partial \Psi}{\partial y}}
    $,
    $
      V
      =
      {\displaystyle \frac{\partial \Phi}{\partial y}}
      =
      - {\displaystyle \frac{\partial \Phi}{\partial x}}
    $.
  }
  A modification of the CBIEM can convert it to become a solver
  for Neumann (and thence mixed) boundary value problems.
  In the Neumann boundary value problem, $ U_{\nu} $, the
  derivative of $U$ across $ \Gamma $, is known instead of
  $U$. Use the Cauchy--Riemann equations to observe that
  $ U_{\nu} = \pm V_{\tau} $ (the tangential component of
  $V$).  The boundary information can be used to construct an
  approximation to $V$, by integration of $ V_{\tau} $ around
  $ \Gamma $, using a suitable zero point (adding in a constant):
  \begin{eqnarray*}
    V \left( \gamma \left( t \right) \right)
    =
    \int_0^t
      V_{\tau}
      \left( \gamma \left( t' \right) \right)
    dt'.
  \end{eqnarray*}
  The same collocation process previously used to approximate
  $V$ can in this case be used to approximate $U$.

  Beyond this, the technique is particularly applicable to free
  boundary problems~\cite{Dold:92}, and may be able to be
  generalised to other elliptic (and possibly other second order)
  operators.
\item
  The method would easily parallelise. The establishment of the
  linear system is computationally expensive, more so for high
  order interpolatory polynomials. This, as well as solution of
  the linear system, would efficiently (geometrically)
  parallelise.
\end{enumerate}

\vfill

\pagebreak

%%%%%%%%%%%%%%%%%%%%%%%%%%%%%%%%%%%%%%%%%%%%%%%%%%%%%%%%%%%%%%%%%%%%%%%%

\appendix

\section{Listing of \protect{\textsc{matlab}} ``.m'' files}
\label{app:code}

\subsection{cbiem.m}
\label{app:cbiem}

\footnotesize
\begin{verbatim}

function Vnnorm = cbiem(CCase, D, sigma, O, alpha);

%function Vnnorm = cbiem(CCase, D, sigma, O, alpha);
%
%       Perform the CBIEM on the Dirichlet problem.
%
%       David  De Wit    March 30  1992  -  January 14  1993

%%%%%%%%%%%%%%%%%%%%%%%%%%%%%%%%%%%%%%%%%%%%%%%%%%%%%%%%%%%%%%%%%%%%%

% 0.0:  Input and other parameters, togther with definitions.

if ~exist('alpha'), alpha = 2/3;  end
if ~exist('O'),     O = 12;        end
if ~exist('sigma'), sigma = 0.32; end
if ~exist('D'),     D = 7;       end
if ~exist('CCase'), CCase = 4;    end

format short e;                                 format compact;

%%%%%%%%%%%%%%%%%%%%%%%%%%%%%%%%%%%%%%%%%%%%%%%%%%%%%%%%%%%%%%%%%%%%%
%
% CCase:Index number to the contour being used:
%       1: Unit circle, with 4 equally spaced artificial corners.
%       2: Chandler's Teardrop, one corner, right angled, at the origin.
%               (This has been reversed to make it ACW. Now different
%               from both GAC and DDW thesis.)
%       3: Kress' ACW Teardrop, one corner, 2 pi/3 - angled.
%       4: Kress' Reentrant contour, 3 pi/2 - angled. Reversed to
%               avoid the branch cut on the negative real axis, and
%               make it ACW in orientation.
%       5: ACW Cardioid. Reentrant contour with 2 pi interior angle.
%       6: ACW Heart. Reentrant contour with 2 pi and 0 interior angles.
%       7: ACW Controlled Cardioid, using trig parameterisation.
%       8: ACW Controlled Cardioid, using polynomial parameterisation.
%       9: Modified Boomerang, with a 5 degree external angle.
%
% D:    Density of the geometric mesh, a positive integer. Choice
%       of D forces N, the size of the linear system being solved,
%       to be N = NC.(D+1)^2.
%
% sigma:Mesh parameter. Ratio of distances of consecutive mesh points
%       from the nearest corner. 0 < sigma < 0.5. Try sigma = 0.25
%       as a starting guess.
%
% O:    Order of the interpolatory polynomial used to aproximate
%       V at the collocation points. Actually the (even) number of
%       nearest points used, thus 2 is linear, 4 is cubic. Must be
%       kept 2 <= O <= NS+1; in practice keep O <= 20. This limits
%       O <= 8 for D = 2, 16 for D = 3, etc.
%
% alpha:Exponent of the true solution of test problem, W = z^{alpha}.
%
% NC:   Numbers of corners/segments of the contour.
%
% NS:   Number of points on each side of the contour. N node points
%       create NS = (N+1)^2 + 1 mesh points on each segment.
%
%%%%%%%%%%%%%%%%%%%%%%%%%%%%%%%%%%%%%%%%%%%%%%%%%%%%%%%%%%%%%%%%%%%%%

% 1.1:  Generate quadrature rule for a geometrically graded mesh, on
%       one segment in the t domain using hprmesh, then insert this
%       into the grid with corners, using rmesh.

numcorn = [4 1 1 1 1 1 1 1 1 1];        NC = numcorn(CCase);

% graded mesh
G = [0 sigma.^(D:-1:1)]';               G = [G; 1-G(D+1:-1:1)];
S = [2:D+2, D+1:-1:2]';                 c = [0:NC]'/NC;
[t, w] = hprmesh(G, S, 0);              [tn, wn] = rmesh(c, t, w, 1);
N = length(tn);

% Uniform mesh
%%% N = (D+1)^2;        tn = [1:N]'/N;  wn = ones(size(tn))/N;

%%%%%%%%%%%%%%%%%%%%%%%%%%%%%%%%%%%%%%%%%%%%%%%%%%%%%%%%%%%%%%%%%%%%%

% 1.2:  Compute the (complex) values of zn and zc (z at the node and
%       collocation points specified by tn and tc, respectively).
%       zc are found as the midpoints of zn in an arc-length sense, by
%       mapping the midpoints of tc to the contour. Also calculate
%       gdot, which is used to modify wn.

tc = tn - diff([0; tn])/2;
ptn = pi*tn;                            ptc = pi*tc;
if (CCase == 1)
        zn = exp(2*i*ptn);              zc = exp(2*i*ptc);
        gdot = 2*pi*i*exp(2*i*ptn);
elseif (CCase == 2) % Chandler's Teardrop
        zn = 2*sin(ptn) - i*sin(2*ptn);
        zc = 2*sin(ptc) - i*sin(2*ptc);
        gdot = 2*pi*(cos(ptn) - i*cos(2*ptn));
        gdot(N) = -2*pi*i;
elseif (CCase == 3) % Kress's Teardrop
        zn = sin(ptn)*2/sqrt(3) - i * sin(2*ptn);
        zc = sin(ptc)*2/sqrt(3) - i * sin(2*ptc);
        gdot = 2*pi*( cos(ptn)/sqrt(3) - i*cos(2*ptn) );
        gdot(N) = -2*pi*i;
elseif (CCase == 4) % Kress's Boomerang
        a = 2/3;
        zn = - a * sin(3*ptn) - i * sin(2*ptn);
        zc = - a * sin(3*ptc) - i * sin(2*ptc);
        gdot = - 2*pi*( (3*a/2)*cos(3*ptn) + i*cos(2*ptn) );
        gdot(N) = -2*pi*i;
elseif (CCase == 5) % Plain Cardioid
        zn = (-1 + cos(2*ptn)).*exp(i*2*ptn);
        zc = (-1 + cos(2*ptc)).*exp(i*2*ptc);
        gdot = -2*pi*exp(i*2*ptn).*(sin(2*ptn) + i*(1-cos(2*ptn)));
        gdot(N) = 0;
elseif (CCase == 6) % Pointed Heart - Silly
        zn = - sin(3*ptn) - i*(sin(2*ptn)).^3;
        zc = - sin(3*ptc) - i*(sin(2*ptc)).^3;
        gdot = - 3*pi*( cos(3*ptn) + i*2*cos(2*ptn).*(sin(2*ptn)).^2 );
        gdot(N/2) = 0;                          gdot(N) = 0;
elseif (CCase == 7) % My Cardioid
        zn = - sin(3*ptn) - 5 * i* tn .* (1 - tn) .* sin(2*ptn);
        zc = - sin(3*ptc) - 5 * i* tc .* (1 - tc) .* sin(2*ptc);
        gdot = - 3 * pi * cos(3*ptn) - 5 * i * ...
                ( (1-2*tn).*sin(2*ptn) + 2*ptn.*(1-tn).*cos(2*ptn) );
        gdot(N) = 0;




elseif (CCase == 8) % My Stupid Polynomial Cardioid
        a = 7/24;
        zn = tn .* (tn - a) .* (tn - 1 + a) .* (tn - 1) + ...
                i * tn.^2 .* (tn - 1/2) .* (tn - 1).^2;
        zc = tc .* (tc - a) .* (tc - 1 + a) .* (tc - 1) + ...
                i * tc.^2 .* (tc - 1/2) .* (tc - 1).^2;
        gdot = (tn - a) .* (tn - 1 + a) .* (tn - 1) + ...
                tn .* (tn-1+a) .* (tn-1) + tn .* (tn-a) .* (tn-1) + ...
                tn .* (tn-a) .* (tn-1+a) + ...
                i * ( 2 * tn .* (tn - 1/2) .* (tn - 1).^2 + ...
                2 * tn.^2 .* (tn - 1/2) .* (tn - 1) + ...
                tn.^2 .* (tn - 1).^2 );
        gdot(N) = 0;
        zn = 100*zn;    zc = 100*zc;    gdot = 100*gdot;
elseif (CCase == 9) % 5 degree external angle.
        deg = 5;
        a = 2/( 3 * tan(deg*pi/360));
        zn = - a * sin(3*ptn) - i * sin(2*ptn);
        zc = - a * sin(3*ptc) - i * sin(2*ptc);
        gdot = - 2*pi*( (3*a/2)*cos(3*ptn) + i*cos(2*ptn) );
        gdot(N) = -2*pi*i;
elseif (CCase == 10) % 20 degree external angle.
        deg = 20;
        a = 2/( 3 * tan(deg*pi/360));
        zn = - a * sin(3*ptn) - i * sin(2*ptn);
        zc = - a * sin(3*ptc) - i * sin(2*ptc);
        gdot = - 2*pi*( (3*a/2)*cos(3*ptn) + i*cos(2*ptn) );
        gdot(N) = -2*pi*i;
end

%rzn = real(zn);                                izn = imag(zn);
%plot(rzn,izn,'-',rzn,izn,'+');         grid;

%plot(tn(1:10),izn(1:10),'-',tn(1:10),izn(1:10),'+');           grid;

wn = wn.*gdot;

%%%%%%%%%%%%%%%%%%%%%%%%%%%%%%%%%%%%%%%%%%%%%%%%%%%%%%%%%%%%%%%%%%%%%%%%

% 2.1:  Set up the complex order N matrices A and B. Establish F, a
%       matrix of indices to be used in calculating B, using Ft, a
%       submatrix of the pattern of indices for one edge. Also compute
%       L, the matrix of the coefficients of the interpolatory
%       polynomial, using the indices contained in F.

oN = ones(N,1);                                 j = O/2;
A = oN*wn.' ./ ( oN*zn.' - zc*oN');             NS = N/NC;

Ft = zeros(NS, O);                              Ft(1,:) = 0:O-1;
Ft(NS,:) = NS-O+1:NS;
Ft(2:j,:) = ones(j-1,1)*Ft(1,:);
Ft(NS-j+1:NS-1,:) = ones(j-1,1)*Ft(NS,:);
for k = j+1:NS-j, Ft(k,:) = Ft(k-1,:) + 1; end
for k = 1:NC, F((k-1)*NS+1:k*NS,:) = Ft + (k-1)*NS; end
F(1:j,1) = N*ones(j,1);








L = ones(size(F));                      o1 = 1:j;       o2 = j+1:N;
for k = 1:O, for v = 1:O
        if (v ~= k)
                tk = tn(F(o2,k));               tv = tn(F(o2,v));
                L(o2, k) = L(o2, k).* (tc(o2) - tv) ./ (tk - tv);

                tk = tn(F(1:j,k));              tv = tn(F(1:j,v));
                if (k == 1), tk = 0; end;       if (v==1), tv = 0; end;
                L(o1, k) = L(o1, k) .* ( tc(o1) - tv ) ./ ( tk - tv );
        end
end, end

%%%%%%%%%%%%%%%%%%%%%%%%%%%%%%%%%%%%%%%%%%%%%%%%%%%%%%%%%%%%%%%%%%%%%%%%

% 2.2:  Compute the matrix B. J is a shift vector, used to create rows
%       of B from rows of L. Multiply the temporary result through by
%       H, then write out C.

H = sum(A.');                                   t1 = ones(j,1);
t2 = [-t1; [0:NS-O-1]'; (NS-O)*t1];             J = zeros(size(H));
for k = 1:NC, J((k-1)*NS+1:k*NS) = t2 + (k-1)*NS; end
B(1:j,1:O-1) = L(1:j,2:O);                      B(1:j,N) = L(1:j,1);
for k = j+1:N, B(k,J(k)+1:J(k)+O) = L(k,:); end
C = B.*(H'*ones(size(H))) - A;

%%%%%%%%%%%%%%%%%%%%%%%%%%%%%%%%%%%%%%%%%%%%%%%%%%%%%%%%%%%%%%%%%%%%%%%%

% 2.3:  Set up and solve for Vn (computed approximation to V at the
%       node points), the linear system:
%               C Vn = A * Utn - sum(A).*Uc
%       Given C, establish the real, order N-1 matrix, Cstar, then
%       calculate Vn, using Vn(N) = 0.

%sprintf('Solving the linear system of size %g', N)
Un = real(funccb(zn, alpha));           Uc = real(funccb(zc, alpha));
d = diag(H)*Uc - A*Un;                  d = -imag(d(2:N) - d(1:N-1));
C = real(C(2:N,1:N-1) - C(1:N-1,1:N-1));
Vn = zeros(size(Un));                   Vn(1:N-1) = C \ d;

%%%%%%%%%%%%%%%%%%%%%%%%%%%%%%%%%%%%%%%%%%%%%%%%%%%%%%%%%%%%%%%%%%%%%%%%

% 2.4:  Error norms. Vne and Vce are the differences between the true
%       and the computed values of V at the node and collocation points
%       respectively. r is the discretisation error in the CBIE. p is
%       the residual in the (above) computation for Vtn, when the true
%       soln at the node points is substituted into the equations.

Vne = imag(funccb(zn, alpha)) - Vn;
Vnnorm = sqrt(abs(wn)'*(Vne.^2));

% Compute the approximation at the 4 points of Kress:

Wn = Un + i * Vn;
z1 = [ 0.1+0*i 0.2+0*i 0.3+0*i 0+0.2*i ];
for j = 1:4
        WK(j) = sum(Wn .* wn ./ (zn - z1(j))) / sum(wn./(zn - z1(j)));
end
WKt = funccb(z1, alpha);                Ek = [N O abs(real(WK - WKt))]

\end{verbatim}

\normalsize

\pagebreak

%%%%%%%%%%%%%%%%%%%%%%%%%%%%%%%%%%%%%%%%%%%%%%%%%%%%%%%%%%%%%%%%%%%%%%%%

\subsection{cint.m}
\label{app:cint}

\footnotesize
\begin{verbatim}

function [rho, C] = cint(sigma);

% function [rho, C] = cint(sigma);
%
%    Contour integration with a geometric h-p grid.
%
%    David  De Wit    July 14  1992  -  July 17  1992

%%%%%%%%%%%%%%%%%%%%%%%%%%%%%%%%%%%%%%%%%%%%%%%%%%%%%%%%%%%%%%%%%%%%%%%%

if ~exist('QCase'), QCase = 2;    end
if ~exist('sigma'), sigma = 0.15; end

DMin = 8;                                       DMax = 15;
format short e;                                 format compact;
QCase = 2;

for D = DMin:DMax
        G = [0 sigma.^(D:-1:1)]';               G = [G; 1-G(D+1:-1:1)];
        S = [2:D+2, D+1:-1:2]';                 NC = 2;
        [t, w] = hprmesh(G, S, QCase, 0);       c = [0:NC]'/NC;
        [tn, wn] = rmesh(c, t, w, 1);           zn = exp(2*pi*i*tn);
        iN(D-DMin+1) = length(zn);
        gdot = 2*pi*i*zn;                       wn = wn.*gdot;
        ier(D-DMin+1) = 1-wn.'*(((zn-1)/i).^(1/2)./zn)/(2*pi*i*sqrt(i));
end
lier = log10(abs(ier))';                        sqiN = sqrt(iN)';

plot(sqiN,lier,'-g',sqiN,lier,'+r');            grid;
title('Contour integration on a h-p geometric grid');
xlabel('sqrt(N)');                              ylabel('log10(Error)');

\end{verbatim}

\normalsize

\pagebreak

%%%%%%%%%%%%%%%%%%%%%%%%%%%%%%%%%%%%%%%%%%%%%%%%%%%%%%%%%%%%%%%%%%%%%%%%

\subsection{funccb.m}
\label{app:funccb}

\footnotesize
\begin{verbatim}

function W = funccb(z, alpha);

% function W = funccb(z, alpha);
%
%       True solution to the Dirichlet problem solved by cbiem.
%
%    David  De Wit    April 13  1992  -  September 27  1992

%%%%%%%%%%%%%%%%%%%%%%%%%%%%%%%%%%%%%%%%%%%%%%%%%%%%%%%%%%%%%%%%%%%%%%%%

W = z.^(alpha);

\end{verbatim}

\normalsize

%%%%%%%%%%%%%%%%%%%%%%%%%%%%%%%%%%%%%%%%%%%%%%%%%%%%%%%%%%%%%%%%%%%%%%%%

\subsection{funcci.m}
\label{app:funcci}

\footnotesize
\begin{verbatim}

function W = funcci(z);

% function W = funcci(z);
%
%       Integrand of the problem solved by cint.
%
%    David  De Wit    July 9  1992  -  July 17  1992

%%%%%%%%%%%%%%%%%%%%%%%%%%%%%%%%%%%%%%%%%%%%%%%%%%%%%%%%%%%%%%%%%%%%%%%%

W = ((z-1)/i).^(1/2)./z;

\end{verbatim}

\normalsize

%%%%%%%%%%%%%%%%%%%%%%%%%%%%%%%%%%%%%%%%%%%%%%%%%%%%%%%%%%%%%%%%%%%%%%%%

\subsection{funchp.m}
\label{app:funchp}

\footnotesize
\begin{verbatim}

function f = funchp(x)

% function f = funchp(x);
%
%       Function being integrated by hpmeth.
%
%    David  De Wit    July 9  1992  -  July 17  1992

%%%%%%%%%%%%%%%%%%%%%%%%%%%%%%%%%%%%%%%%%%%%%%%%%%%%%%%%%%%%%%%%%%%%%%%%

f = 1 - 3/2*sqrt(x);

\end{verbatim}

\normalsize

\pagebreak

%%%%%%%%%%%%%%%%%%%%%%%%%%%%%%%%%%%%%%%%%%%%%%%%%%%%%%%%%%%%%%%%%%%%%%%%

\subsection{gettw.m}
\label{app:gettw}

\footnotesize
\begin{verbatim}

function [QRt, QRw] = gettw(R, QCase);

% function [QRt, QRw] = gettw(R, QCase);
%
%       Get tables of nodes and weights for quadrature rules. User
%       inputs maximum number of points required, and the type required.
%       The default type is Gauss--Lobatto (QCase = 2), as it is of
%       higher order than Newton--Cotes.
%
%       David  De Wit    July 13  1992  -  July 17  1992

%%%%%%%%%%%%%%%%%%%%%%%%%%%%%%%%%%%%%%%%%%%%%%%%%%%%%%%%%%%%%%%%%%%%%%%%

if ~exist('R'),     R = 20;    end
if ~exist('QCase'), QCase = 2; end

%%%%%%%%%%%%%%%%%%%%%%%%%%%%%%%%%%%%%%%%%%%%%%%%%%%%%%%%%%%%%%%%%%%%%%%%

if (QCase == 1)
        if (R > 10)
                R = 10;
                sprintf('R too large for Newton--Cotes. Now R = 10.\n');
        end
        QRw = [
                1  1  1   7  19   41   751    989   2857    16067;
                1  4  3  32  75  216  3577   5888  15741   106300;
                0  1  3  12  50   27  1323   -928   1080   -48525;
                0  0  1  32  50  272  2989  10496  19344   272400;
                0  0  0   7  75   27  2989  -4540   5778  -260550;
                0  0  0   0  19  216  1323  10496   5778   427368;
                0  0  0   0   0   41  3577   -928  19344  -260550;
                0  0  0   0   0    0   751   5888   1080   272400;
                0  0  0   0   0    0     0    989  15741   -48525;
                0  0  0   0   0    0     0      0   2857   106300;
                0  0  0   0   0    0     0      0      0    16067
        ]
        QRw = QRw(1:R+1,1:R);                   QRt = zeros(QRw);
        for j = 1:R
                QRw(1:j+1,j) = QRw(1:j+1,j)/sum(QRw(:,j));
                QRt(1:j+1,j) = [0:j]'/j;
        end
elseif (QCase == 2)
        QRt = zeros(R+1,R);                     QRw = QRt;
        for j = 2:R+1
                [QRt(1:j,j-1), QRw(1:j,j-1)] = lobatto(j,0,1);
        end
end

\end{verbatim}

\normalsize

\pagebreak

%%%%%%%%%%%%%%%%%%%%%%%%%%%%%%%%%%%%%%%%%%%%%%%%%%%%%%%%%%%%%%%%%%%%%%%%

\subsection{hpmeth.m}
\label{app:hpmeth}

\footnotesize
\begin{verbatim}

function [lEip] = hpmeth(DMax, GMax, p, sigma, QCase);

% function [lEip] = hpmeth(DMax, GMax, p, sigma, QCase);
%
%       Experiment with h-p integration methods.
%
%       David  De Wit    July 9  1992  -  July 17  1992

%%%%%%%%%%%%%%%%%%%%%%%%%%%%%%%%%%%%%%%%%%%%%%%%%%%%%%%%%%%%%%%%%%%%%%%%

if ~exist('QCase'), QCase = 2;    end
if ~exist('sigma'), sigma = 0.15; end
if ~exist('p'),     p = 6;        end
if ~exist('GMax'),  GMax = 6;     end
if ~exist('DMax'),  DMax = 19;    end

%%%%%%%%%%%%%%%%%%%%%%%%%%%%%%%%%%%%%%%%%%%%%%%%%%%%%%%%%%%%%%%%%%%%%%%%

% Geometrically-graded h-p method.

Ehp = zeros(DMax,1);                            vN = Ehp;
for D = 1:DMax
        G = [0 sigma.^(D:-1:0)]';               S = [2:D+2]';
        [t, w] = hprmesh(G, S, QCase, 0);
        Ehp(D) = funchp(t)'*w;                  vNhp(D) = length(t);
end
lvNhp = log10(vNhp);                            lEhp = log10(Ehp);

%%%%%%%%%%%%%%%%%%%%%%%%%%%%%%%%%%%%%%%%%%%%%%%%%%%%%%%%%%%%%%%%%%%%%%%%

% Various linear h and p methods. Variable g, constant p.

Eip = zeros(DMax,GMax);
for D = 1:DMax
        N = 2*D;
        for g = 1:GMax
                G = ([0:N]'/N).^g;
                [t, w] = hprmesh(G, p, QCase, 0);
                vN(D) = length(t);              Eip(D,g) = funchp(t)'*w;
        end
end

lvN = log10(vN);                                lEip = log10(Eip);

sprintf('Order %f, slopes of the blue lines are approximately', p)
(lEip(DMax-1,:) - lEip(DMax,:))./(lvN(DMax-1)-lvN(DMax))

plot(lvNhp,lEhp,'+r',lvNhp,lEhp,'-g',lvN,lEip,'+r',lvN,lEip,'-b');
xlabel('log10(points)');                ylabel('log10(error)');
text(0.9,0.7,'g = 1','sc');             text(0.9,0.57,'g = 2','sc');
text(0.9,0.45,'g = 3','sc');            text(0.9,0.35,'g = 4','sc');
text(0.9,0.27,'g = 5','sc');            text(0.9,0.2,'g = 6','sc');
text(0.7,0.15,'h-p method','sc');       grid;

%%%%%%%%%%%%%%%%%%%%%%%%%%%%%%%%%%%%%%%%%%%%%%%%%%%%%%%%%%%%%%%%%%%%%%%%






% Compare h methods for quadrature rules of various p (# points).

OMax = 12;
for p = 2:2:OMax
        g = p;                                  j = p/2 + 1;
        for D = 1:DMax
                N = 2*D;                        G = ([0:N]'/N).^g;
                [t, w] = hprmesh(G, j, QCase, 0);
                Ehhp(D,j-1) = funchp(t)'*w;     hhN(D) = length(t);
        end
end

lhhN = log10(hhN);                              lEhhp = log10(Ehhp);

sprintf('Slopes of the blue lines are approximately')
(lEhhp(DMax-1,:) - lEhhp(DMax,:))./(lhhN(DMax-1)-lhhN(DMax))

plot(lvNhp,lEhp,'+r',lvNhp,lEhp,'-g',lhhN,lEhhp,'+r',lhhN,lEhhp,'-b');
xlabel('log10(points)');                ylabel('log10(error)');
text(0.9,0.76,'p = 1','sc');            text(0.9,0.62,'p = 3','sc');
text(0.9,0.5,'p = 5','sc');             text(0.9,0.4,'p = 7','sc');
text(0.9,0.32,'p = 9','sc');            text(0.85,0.2,'p = 11','sc');
text(0.7,0.15,'h-p method','sc');       grid;

\end{verbatim}

\normalsize

\pagebreak

%%%%%%%%%%%%%%%%%%%%%%%%%%%%%%%%%%%%%%%%%%%%%%%%%%%%%%%%%%%%%%%%%%%%%%%%

\subsection{hprmesh.m}
\label{app:hprmesh}

\footnotesize
\begin{verbatim}

function [t, w] = hprmesh(G, S, IClosed)

% function [t, w] = hprmesh(G, S, IClosed)
%
%       Create a new quadrature rule, based on a mesh G, where between
%       points G(i) and G(i+1) is a (closed) quadrature rule of
%       Gauss--Lobatto type on S(i) points, including the 2 end points
%       G(i) and G(i+1); for i = 1:length(G)-1.  If IClosed is 1, then
%       the contour is closed, and the ends are tied together.  This
%       function is a generalisation of rmesh.
%
%       David  De Wit   July 13  1992  -  December 2  1992

%%%%%%%%%%%%%%%%%%%%%%%%%%%%%%%%%%%%%%%%%%%%%%%%%%%%%%%%%%%%%%%%%%%%%%%%

if ~exist('G'),       sigma = 0.1; G = [0 sigma.^(3:-1:1) 1]'; end
if ~exist('S'),       S = [2:5]';                              end
if ~exist('IClosed'), IClosed = 1;                             end

lG = length(G);                         lS = length(S);
if ((lG ~= lS+1) & (lS ~= 1))
        sprintf('hprmesh: Danger l(G) = %f, l(S) = %f', lG, lS)
end

% Set up S for rules with a constant integration rule.
if (lS == 1), S = ones(lG-1,1)*S; end

dG = diff(G);                                   S = S - 1;
N = max(S);                                     R = length(dG);

% Obtain the nodes and weights in a table
QRt = zeros(N+1,N);                     QRw = QRt;
for j = 2:N+1
        [QRt(1:j,j-1), QRw(1:j,j-1)] = lobatto(j,0,1);
end

% Play with the table
QRw(N+1,:) = diag(QRw(2:N+1,1:N))';             QRt(N+1,:) = ones(1,N);
for i = 2:N, for j = 1:i-1, QRt(i,j) = NaN; QRw(i,j) = NaN; end, end
tt = QRt(:,S);                                   tw = QRw(:,S);
j = ones(N+1,1);                                 tw = tw.*(j*dG');
tt = tt.*(j*dG') + j*G(1:R)';
tw(1,2:R) = tw(1,2:R) + tw(N+1,1:R-1);
tw1 = tw(1:N,:);                                tt1 = tt(1:N,:);
t = [tt1(:); tt(N+1,R)];                        t(isnan(t)) = [];
w = [tw1(:); tw(N+1,R)];                        w(isnan(w)) = [];

if (IClosed == 1)
        N = length(t);
        w = [w(2:N-1); w(1)+w(N)];              t = t(2:N);
end

\end{verbatim}

\normalsize

\pagebreak

%%%%%%%%%%%%%%%%%%%%%%%%%%%%%%%%%%%%%%%%%%%%%%%%%%%%%%%%%%%%%%%%%%%%%%%%

\subsection{lobatto.m}
\label{app:lobatto}

\footnotesize
\begin{verbatim}

function [x, w] = lobatto(n, a, b)

% function [x, w] = lobatto(n, a, b)
%
%       Return the weights w and points x of the n-point Gauss--Lobatto
%       quadrature rule on the interval [a, b].
%       See G. H. Golub, SIAM Review 1973 p 318.
%
%       Graeme Chandler         July 1992

%%%%%%%%%%%%%%%%%%%%%%%%%%%%%%%%%%%%%%%%%%%%%%%%%%%%%%%%%%%%%%%%%%%%%%%%

n = round(n);

if (n == 2)
        x = [a; b];                             w = [1; 1]*(b-a)/2;
elseif (n == 3)
        x = [a; a+(b-a)/2; b];                  w = [1; 4; 1]*(b-a)/6;
elseif (n >= 4)
        nn = n-1;                               m = 1:2:2*nn-1;
        m = (1:nn-1) ./ sqrt(m(1:nn-1) .* m(2:nn));
        J = (diag(m,-1)+diag(m,1));
        I = eye(nn);                            en = (1:nn)' == nn;
        gam = (J + I)\en;                       mu = (J - I)\en;
        sol = [1 -gam(nn); 1 -mu(nn)]\[-1; 1];
        alpha = sol(1);                         beta = sqrt(sol(2));
        [ww,xx] = eig([J beta*en; beta*en' alpha]);
        [xx, i] = sort(diag(xx));
        w = ww(1,i)'.^2 * (b-a);
        x = [a; (a+b)/2+(b-a)*xx(2:nn)/2; b];
end

\end{verbatim}

\normalsize

\pagebreak

%%%%%%%%%%%%%%%%%%%%%%%%%%%%%%%%%%%%%%%%%%%%%%%%%%%%%%%%%%%%%%%%%%%%%%%%

\subsection{rmesh.m}
\label{app:rmesh}

\footnotesize
\begin{verbatim}

function [t, w] = rmesh(G, t1, w1, IClosed);

% function [t, w] = rmesh(G, t1, w1, IClosed);
%
%       Create a new quadrature rule, based on a mesh G, where a
%       (closed) quadrature rule (t1, w1) is inserted over each
%       interval of G.  If IClosed is 1, then the contour is closed,
%       and the ends are tied together.  This function is generalised
%       into hprmesh. Originally conceived by Graeme Chandler.
%
%       David De Wit    July 13  1992  -  September 6  1992
%
%%%%%%%%%%%%%%%%%%%%%%%%%%%%%%%%%%%%%%%%%%%%%%%%%%%%%%%%%%%%%%%%%%%%%%%%

if (length(G) < 2), return; end

if (length(G) ~= 2)
        h = diff(G);                            tw = w1*h';
        n = length(h);                          m = length(w1);
        tw(1,2:n) = tw(1,2:n) + tw(m,1:n-1);    tw1 = tw(1:m-1,:);
        w = [tw1(:); tw(m,n)];
        tt = t1(1:m-1)*h' + ones(m-1,1)*G(1:n)';
        t = [tt(:); G(n+1)];
else
        t = t1;                                 w = w1;
end

if (IClosed == 1)
        N = length(t);
        w = [w(2:N-1); w(1)+w(N)];              t = t(2:N);
end

\end{verbatim}

\normalsize

\pagebreak

%%%%%%%%%%%%%%%%%%%%%%%%%%%%%%%%%%%%%%%%%%%%%%%%%%%%%%%%%%%%%%%%%%%%%%%%

\subsection{testcb.m}
\label{app:testcb}

\footnotesize
\begin{verbatim}

function N = testcb(C, sigma, Dmin, Dmax, Omin, Omax, alpha)

% function N = testcb(C, sigma, Dmin, Dmax, Omin, Omax, alpha)
%
%       Run cbiem for various parameters, and tabulate results.
%
%       David  De Wit    May 12  1992  -  December 21  1992

%%%%%%%%%%%%%%%%%%%%%%%%%%%%%%%%%%%%%%%%%%%%%%%%%%%%%%%%%%%%%%%%%%%%%%%%

if ~exist('alpha'), alpha = 1/2; end
if ~exist('Omax'), Omax = 16; end
if ~exist('Omin'), Omin = 8; end
if ~exist('Dmax'), Dmax = 19; end
if ~exist('Dmin'), Dmin = 16; end
if ~exist('sigma'), sigma = 0.32; end
if ~exist('C'), C = 7; end
format short e;                                  format compact

%%%%%%%%%%%%%%%%%%%%%%%%%%%%%%%%%%%%%%%%%%%%%%%%%%%%%%%%%%%%%%%%%%%%%%%%

for D = Dmin:Dmax
        for O = Omin:2:Omax
                P = cbiem(C, D, sigma, O, alpha);
                if ((P ~= Inf) & (P ~= NaN))
                        N(D-Dmin+1, (O-Omin)/2+1) = P;
                else
                        N(D-Dmin+1:Dmax-Dmin+1,:) = ...
                                Inf*ones(Dmax-D+1,(Omax-O)/2 + 1);
                        return
                end
        end
        if (min(N(D-Dmin+1,:)) > 1)
                N(D-Dmin+2:Dmax-Dmin+1,:) = ...
                        Inf*ones(Dmax-D,(Omax-Omin)/2+1);
                return
        end
end

\end{verbatim}

\normalsize

%%%%%%%%%%%%%%%%%%%%%%%%%%%%%%%%%%%%%%%%%%%%%%%%%%%%%%%%%%%%%%%%%%%%%%%%

\pagebreak

%%%%%%%%%%%%%%%%%%%%%%%%%%%%%%%%%%%%%%%%%%%%%%%%%%%%%%%%%%%%%%%%%%%%%%%%

\bibliographystyle{plain}
\bibliography{DeWit92}

\end{document}